\documentclass[11pt]{amsart}
\usepackage{amssymb}
\usepackage{amscd}
\usepackage[all]{xy}
\usepackage{amsthm}

\newtheorem{Thm}{Theorem}[section]
\newtheorem{Lem}[Thm]{Lemma}
\newtheorem{Cor}[Thm]{Corollary}
\newtheorem{Prop}[Thm]{Proposition}
\newtheorem{Conj}[Thm]{Conjecture}

\newcommand{\A}{\mathbb{A}}
\newcommand{\D}{\mathbb{D}}
\newcommand{\E}{\mathbb{E}}

\newcommand{\Z}{\mathbb{Z}}

\newcommand{\N}{\mathbb{N}}
\newcommand{\C}{\mathbb{C}}

\newcommand{\Lam}{\Lambda}

\newcommand{\id}{1\kern -.35em 1}

\newcommand{\G}{{\mathcal G}}
\newcommand{\g}{\mathfrak{g}}
\newcommand{\SB}{{\mathcal S}}
\newcommand{\CB}{{\mathcal C}}
\newcommand{\MB}{{\mathcal M}}
\newcommand{\T}{{\mathcal T}}

\newcommand{\orb}{\mathcal O}
\newcommand{\md}{\operatorname{mod}}

\newcommand{\Hom}{\operatorname{Hom}}
\newcommand{\Ext}{\operatorname{Ext}}
\newcommand{\End}{\operatorname{End}}
\newcommand{\dm}{{\rm dim}\,}
\newcommand{\GL}{{\rm G}}
\newcommand{\field}{K}
\newcommand{\rep}{\operatorname{rep}}
\newcommand{\Ima}{\operatorname{Im}}
\newcommand{\Ker}{\operatorname{Ker}}
\newcommand{\Coker}{\operatorname{Coker}}

\newcommand{\la}{\Lambda}

\newcommand{\bil}[1]{\langle #1\rangle}

\newcommand{\dimv}{\underline{\dim}}

\newcommand{\op}{{\rm op}}

\newcommand{\gldim}{\operatorname{gl.dim}}
\newcommand{\pdim}{\operatorname{proj.dim}}

\newcommand{\repdim}{\operatorname{rep.dim}}
\newcommand{\domdim}{\operatorname{dom.dim}}

\newcommand{\n}{\mathfrak{n}}
\newcommand{\be}{\beta}

\newcommand{\irr}{\operatorname{Irr}}
\newcommand{\xb}{\mathbf{x}}
\newcommand{\ib}{\mathbf{i}}
\newcommand{\ab}{\mathbf{a}}
\newcommand{\tB}{\widetilde{B}}

\begin{document}

\bigskip
\title{Rigid modules over preprojective algebras}

\author{Christof Gei{\ss}}
\address{Christof Gei{ss}\newline
Instituto de Matem\'aticas, UNAM\newline
Ciudad Universitaria\newline
04510 Mexico D.F.\newline
Mexico}
\email{christof@math.unam.mx}

\author{Bernard Leclerc}
\address{Bernard Leclerc\newline
Laboratoire LMNO\newline
Universit\'e de Caen\newline
F-14032 Caen Cedex\newline
France}
\email{leclerc@math.unicaen.fr}

\author{Jan Schr\"oer}
\address{Jan Schr\"oer\newline
Mathematisches Institut\newline
Universit\"at Bonn\newline
Beringstr. 1\newline
D-53115 Bonn\newline
Germany}
\email{schroer@math.uni-bonn.de}

\thanks{
J. Schr\"oer was supported by a research fellowship from 
the DFG (Deutsche Forschungsgemeinschaft).
He also thanks the Laboratoire LMNO (Caen) for an invitation
in May 2004 during which this work was started. 
}


\begin{abstract}
Let $\la$ be a preprojective algebra of simply laced Dynkin type
$\Delta$.
We study maximal rigid $\la$-modules, their endomorphism algebras
and a mutation operation on these modules.
This leads to a representation-theoretic construction
of the cluster algebra structure on the ring $\C[N]$ 
of polynomial functions on a maximal unipotent subgroup $N$
of a complex Lie group of type $\Delta$.
As an application we obtain that all cluster monomials of 
$\C[N]$ belong to the dual semicanonical basis.
\end{abstract}

\maketitle

\bigskip
\setcounter{tocdepth}{1}
\tableofcontents



\section{Introduction}


\subsection{}
Preprojective algebras were introduced by Gelfand and Ponomarev
in 1979 \cite{GP}, 
and since then played an important role in representation
theory and Lie theory.
For example, let $U(\n)$ be the enveloping algebra of
a maximal nilpotent subalgebra $\n$ of a simple Lie algebra $\g$
of type $\A, \D, \E$, and let $\Lam$ denote the preprojective
algebra associated to the Dynkin diagram of $\g$.
In \cite{Lu}, \cite{Lu2}, Lusztig has given a geometric construction
of $U(\n)$ as an algebra $(\MB,*)$ of constructible
functions on varieties of finite-dimensional $\Lam$-modules, 
where $*$ is a convolution product inspired by Ringel's
multiplication for Hall algebras 
\cite{Ri2}. 
This yields a new basis $\SB$ of $U(\n)$ given by 
the irreducible components of these varieties of modules,
called the {\it semicanonical basis}.

\subsection{}
Cluster algebras were introduced by Fomin and Zelevinsky
\cite{FZ} to provide, among other things, an algebraic and 
combinatorial framework for the study of canonical bases
and of total positivity.
Of particular interest is the algebra $\C[N]$ of polynomial
functions on a unipotent group $N$ with Lie algebra $\n$.
One can identify $\C[N]$ with the graded dual $U(\n)^*$ and thus
think of the dual $\SB^*$ of $\SB$   
as a basis  of $\C[N]$. We call $\SB^*$ the {\it dual 
semicanonical basis}.
It follows from \cite{BFZ} that $\C[N]$ has a natural (upper) 
cluster algebra structure. This consists of a distinguished
family of regular functions on $N$, grouped into subsets called
clusters. It is generated inductively from an initial cluster
defined in a combinatorial way. 
In type $\A_n \ (n\le 4)$, there are finitely many clusters
and their elements, the cluster variables, can all be     
described explicitly. 
Moreover the set of all cluster monomials  
(that is, monomials in the cluster variables supported on
a single cluster) coincides with $\SB^*$.
In general, though, there are infinitely many clusters
and cluster variables in $\C[N]$, and very little is known about them.

\subsection{} \label{lift}
To a finite dimensional $\Lam$-module $M$ one can attach 
the linear form $\delta_M\in\MB^*$ which maps a constructible 
function $f\in\MB$ to its evaluation $f(M)$ at $M$.
Under the isomorphism $\MB^* \cong U(\n)^* \cong \C[N]$,
$\delta_M$ gets identified to a regular function 
$\varphi_M\in\C[N]$. 
Thus, to study special elements of $\C[N]$, like cluster monomials,
we may try to lift them to $\md(\Lam)$ via the map $M \mapsto \varphi_M$. 
For example, by construction, the element of $\SB^*$ attached to 
an irreducible 
component $Z$ of a variety of $\Lam$-modules is equal to
$\varphi_M$, where $M$ is a ``generic module'' in $Z$. 

\subsection{}
A $\la$-module $M$ is called {\it rigid} provided $\Ext_\la^1(M,M)=0$.
We characterize rigid $\la$-modules as the modules having an open orbit in the
corresponding module variety.
In particular, the closure of such an orbit is an irreducible component.

In the first part of this paper, 
we show that the endomorphism algebras of rigid modules
have astonishing properties
which we believe are interesting in themselves.
In particular, maximal rigid modules are examples of maximal 1-orthogonal
modules, which play a role in the 
higher dimensional Auslander-Reiten theory recently developed by
Iyama \cite{I}, \cite{I1}. 
This yields a direct link between preprojective algebras and classical 
tilting theory.

In the second part, we use these results to
show that the operation of {\it mutation}
involved in the cluster algebra structure of $\C[N]$ 
can be entirely understood in terms 
of maximal rigid $\la$-modules.
More precisely we define a mutation operation on maximal
rigid modules, and, taking into account the results of \cite{GLS}, 
\cite{GLS2}, \cite{GLS3}, we prove that this gives a lifting
of the cluster structure to the category $\md(\Lam)$.
In particular, this implies that all cluster monomials of $\C[N]$ 
belong to $\SB^*$. 
This theorem establishes for the first time a bridge between
Lusztig's geometric construction of canonical bases and
Fomin and Zelevinsky's approach to this topic via cluster algebras.

\subsection{}
Our way of understanding the cluster algebra $\C[N]$
via the category 
$\md(\Lam)$ is very similar to the approach of  
\cite{BMRRT}, \cite{BMR}, \cite{BMR2}, 
\cite{CC}, \cite{CK}, \cite{K}
which study cluster algebras attached to quivers via some
new {\it cluster categories}.
There are however two main differences.

First, the theory of cluster categories 
relies on the well developed representation 
theory of hereditary algebras, and covers only a special
class of cluster algebras called {\it acyclic}.
With the exception of Dynkin types
$\A_n$ with $n \le 4$, the cluster algebras $\C[N]$ 
are not believed to belong to this class. 
This indicates that hereditary algebras cannot be used in this context.
In contrast, 
we use the preprojective algebras which have infinite global 
dimension and whose representation theory is much less developed.

Secondly, cluster categories are not abelian categories
but only triangulated categories, defined as orbit categories
of derived categories of representations of quivers.  
In our approach, we just use the concrete abelian category 
$\md(\Lam)$.


\section{Main results}


\subsection{}
Let $A$ be a finite-dimensional algebra over an algebraically
closed field $\field$.
By $\md(A)$ we denote the category of finite-dimensional 
left $A$-modules.
If not mentioned otherwise, modules are assumed to be 
left modules.
In this article we only consider finite-dimensional modules.
We often do not distinguish between a module and its isomorphism
class.
For an $A$-module $M$ let ${\rm add}(M)$ be the full subcategory of
$\md(A)$ formed by all modules isomorphic to direct summands of finite direct
sums of copies of $M$.
The opposite algebra of $A$ is denoted by $A^\op$.
Let 
$$
{\rm D} = \Hom_\field(-,\field)\colon \md(A) \to \md(A^\op) 
$$ 
be the usual duality functor.

If $f\colon U \to V$ and $g\colon V \to W$ are morphisms, then the composition
is denoted by $gf\colon U \to W$.
With this convention concerning compositions of homomorphisms we get
that $\Hom_A(M,N)$ is a left $\End_A(N)$-module and a right
$\End_A(M)$-module.

For natural numbers $a \leq b$ let $[a,b] = \{ i \in \N \mid a \le i \le b\}$.

\subsection{}
An $A$-module $T$ is called {\it rigid} 
if $\Ext_A^1(T,T) = 0$.
A rigid module $T$ is {\it maximal} if 
any indecomposable $A$-module $T'$ such that $T \oplus T'$ 
is rigid, is isomorphic to a direct summand of $T$
(in other words, $T' \in {\rm add}(T)$).

\subsection{}
Throughout the article let $Q = (Q_0,Q_1,s,t)$ be a Dynkin quiver of
simply laced type 
\[ 
\Delta \in
\{ \A_n \, (n \ge 2), \, \D_n \, (n \ge 4), \, \E_n \, (n = 6,7,8) \}.
\]
Thus $Q$ is given by a simply laced Dynkin diagram together with
an arbitrary orientation on the edges.
Here $Q_0$ and $Q_1$ denote the set of vertices and arrows of $Q$,
respectively.
For an arrow $\alpha\colon i \to j$ in $Q$ let $s(\alpha) = i$ and
$t(\alpha) = j$ be its starting and terminal vertex, respectively.
By $n$ we always denote the number of vertices of $Q$.
Note that we exclude the trivial case $\A_1$.

Let $\overline{Q}$ be the {\it double quiver} of $Q$, which is obtained from
$Q$ by adding an arrow $\alpha^*\colon j \to i$ whenever there is an 
arrow $\alpha\colon i \to j$ in $Q$.
The {\it preprojective algebra} associated to $Q$ is defined as
\[
\la = \la_Q = \field \overline{Q}/(c)
\]
where $(c)$ is the ideal generated by the element
\[
c = \sum_{\alpha \in Q_1} (\alpha^*\alpha - \alpha\alpha^*),
\]
and  $\field \overline{Q}$ is the path algebra associated to
$\overline{Q}$,
see \cite{Ri3}.
Since $Q$ is a Dynkin quiver it follows that $\la$ is a 
finite-dimensional selfinjective algebra.
One can easily show that $\la$ does not
depend on the orientation of $Q$.
More precisely, if $Q$ and $Q'$ are Dynkin quivers of the same Dynkin
type $\Delta$, then $\Lambda_Q$ and $\Lambda_{Q'}$ are isomorphic algebras.

Let $r$ be the set of positive roots of $Q$, or equivalently, let
$r$ be the number of isomorphism classes of indecomposable
representations of $Q$, compare \cite{Ga}.
Here are all possible values for $r$:
\begin{center}
{\renewcommand{\arraystretch}{1.5}
\begin{tabular}{c||c|c|c|c|c}
$Q$  & $\A_n$             & $\D_n$  & $\E_6$ & $\E_7$ & $\E_8$ \\ \hline
$r$  & $\frac{n(n+1)}{2}$ & $n^2-n$ & $36$   & $63$   & $120$
\end{tabular}
}
\end{center}
\medskip\noindent
For a module $M$ let $\Sigma(M)$ be the number of 
isomorphism classes of indecomposable direct summands of 
$M$.
The following theorem is proved in \cite{GSc}:

\begin{Thm}\label{upperbound}
For any rigid $\la$-module $T$ we have $\Sigma(T) \leq r$.
\end{Thm}

We call a rigid $\la$-module $T$ {\it complete} if $\Sigma(T) = r$.
It follows from the definitions and from Theorem \ref{upperbound} that
any complete rigid module is also maximal rigid. 
In Theorem \ref{introthm1} we will show the converse, namely that
$\Sigma(T) = r$ for any maximal rigid module $T$ .

\subsection{}
Let $A$ be a finite-dimensional algebra.
Following Iyama we call an 
additive full subcategory ${\mathcal T}$ of $\md(A)$  
{\it maximal 1-orthogonal} if for every $A$-module $M$ 
the following are equivalent:
\begin{itemize}

\item 
$M \in {\mathcal T}$;

\item
$\Ext_A^1(M,T) = 0$ for all $T \in {\mathcal T}$;

\item 
$\Ext_A^1(T,M) = 0$ for all $T \in {\mathcal T}$.

\end{itemize}
An $A$-module $T$ is called {\it maximal 1-orthogonal} if ${\rm add}(T)$
is maximal 1-orthogonal.

An $A$-module $C$ is a
{\it generator} (resp. {\it cogenerator}) of $\md(A)$ if for every 
$A$-module $M$ 
there exists some $m \ge 1$ and an epimorphism $C^m \to M$ (resp.
a monomorphism $M \to C^m$).
One calls $C$ a {\it generator-cogenerator} if it is both a generator and
a cogenerator.
It follows that $C$ is a generator (resp. cogenerator) if and only if
all indecomposable projective (resp. injective) $A$-modules
occur as direct summands of $C$, up to isomorphism.

It follows from the definitions that any maximal 1-orthogonal 
$A$-module $T$ is a 
generator-cogenerator of $\md(A)$.
It also follows that $T$ is rigid.
In general, maximal 1-orthogonal modules need not exist.

The {\it global dimension} $\gldim(A)$ of an algebra $A$ is the supremum 
over all projective dimensions $\pdim(M)$ of all $A$-modules $M$ in 
$\md(A)$.
The {\it representation dimension} of $A$ is defined as
\[
\repdim(A) = \inf \{ \gldim(\End_A(C)) \mid 
C \text{ a generator-cogenerator of } \md(A) \}.
\]
Auslander proved that $A$ is of finite representation type if and only if
$\repdim(A) \le 2$, see \cite[p.559]{A}. 
Iyama \cite{I0} showed that $\repdim(A)$ is always finite.
Next, let
\[
0 \to A \to I_0 \to I_1 \to I_2 \to \cdots
\]
be a minimal injective resolution of $A$.
Then the {\it dominant dimension} of $A$ is defined as
\[
\domdim(A) = \inf \{ i \ge 0 \mid \text{$I_i$ is non-projective} \}.
\]

Let $M$ be a finite-dimensional $A$-module.
Thus $M = M_1^{n_1} \oplus \cdots \oplus M_t^{n_t}$
where the $M_i$ are pairwise non-isomorphic indecomposable
modules and $n_i = [M:M_i] \ge 1$ is the 
multiplicity of $M_i$ in $M$.
The endomorphism algebra $\End_A(M)$ is Morita equivalent to an algebra 
$\field \Gamma_M/I$, where
$\Gamma_M$ is some uniquely determined
finite quiver and $I$ is some admissible ideal in the path algebra
$\field \Gamma_M$.
(An ideal $I$ in a path algebra $\field \Gamma$ is called {\it admissible}
if $I$ is generated by a set of elements of the form
$\sum_{i=1}^m \lambda_i p_i$, where the $p_i$ are paths of length at
least two in $\Gamma$ and the $\lambda_i$ are in $\field$, 
and if $\field \Gamma/I$
is finite-dimensional.)

The module $M$ is called {\it basic} if $n_i = 1$ for all $i$.
In fact we have 
$\End_A(M_1 \oplus \cdots \oplus M_t) \cong \field \Gamma_M/I$.
Thus, if $M'$ is the basic module
$M_1 \oplus \cdots \oplus M_t$, then the categories
$\md(\End_A(M))$ and $\md(\End_A(M'))$ are equivalent, and we can 
restrict to the study of $\md(\End_A(M'))$.

We call $\Gamma_M$ the {\it quiver} of $\End_A(M)$.
The vertices $1, \ldots, t$ of $\Gamma_M$ correspond to the modules $M_i$.
By $S_{M_i}$ or just $S_i$ we denote the simple 
$\End_A(M)$-module corresponding to $i$.
The indecomposable projective $\End_A(M)$-module $P_i$ with top
$S_i$ is just $\Hom_A(M_i,M)$.

Our first main result shows that endomorphism algebras of maximal
rigid modules over preprojective algebras have surprisingly
nice properties:

\begin{Thm}\label{introthm1}
Let $\la$ be a preprojective algebra of Dynkin type
$\Delta$.
For a $\la$-module $T$ the following are equivalent:
\begin{itemize}

\item
$T$ is maximal rigid;

\item
$T$ is complete rigid;

\item
$T$ is maximal 1-orthogonal.

\end{itemize}
If $T$ satisfies one of the above equivalent conditions, then
the following hold:
\begin{itemize}

\item
$\gldim(\End_\la(T)) = 3$;

\item
$\domdim(\End_\la(T)) = 3$;

\item
The quiver $\Gamma_T$ of $\End_\la(T)$ has no sinks, no sources, no loops 
and no 2-cycles.

\end{itemize}
\end{Thm}

\begin{Cor}\label{introthm1cor1}
$\repdim(\la) \le 3$.
\end{Cor}

The following is a consequence of Theorem \ref{introthm1}, see
Proposition \ref{homfunctor} for a proof:

\begin{Cor}\label{introthm1cor2}
For a maximal rigid $\la$-module $T$
the functor
\[
\Hom_\la(-,T)\colon \md(\la) \to \md(\End_\la(T))
\]
is fully faithful, and its image is the category of $\End_\la(T)$-modules
of projective dimension at most one.
\end{Cor}

The proof of Theorem \ref{introthm1} uses 
the following variation of a result
by Iyama (see Theorem \ref{Iyama3} below).
The definition of a tilting module can be found in Section 
\ref{tiltingsection}.

\begin{Prop}
Let $T_1$ and $T_2$ be maximal rigid $\la$-modules.
Then $T = \Hom_\la(T_2,T_1)$ is a tilting module over $\End_\la(T_1)$,
and we have
\[
\End_{\End_\la(T_1)}(T) \cong \End_\la(T_2)^\op.
\]
In particular, $\End_\la(T_1)$ and $\End_\la(T_2)$ are derived equivalent.
\end{Prop}

\subsection{}\label{matmutation}
If $\widetilde{B}$ is any $r \times (r-n)$-matrix, then
the {\it principal part} $B$ of 
$\widetilde{B}$ is obtained from 
$\widetilde{B}$ by deleting the last $n$ rows.
The following definition is due to Fomin and Zelevinsky \cite{FZ}:
Given some $k \in [1,r-n]$ define a new $r \times (r-n)$-matrix 
$\mu_k(\widetilde{B}) = (b_{ij}')$ by
\[
b_{ij}' =
\begin{cases}
-b_{ij} & \text{if $i=k$ or $j=k$},\\
b_{ij} + \dfrac{|b_{ik}|b_{kj} + b_{ik}|b_{kj}|}{2} & \text{otherwise},
\end{cases}
\]
where $i \in [1,r]$ and $j \in [1,r-n]$.
One calls $\mu_k(\widetilde{B})$ a {\it mutation} of $\widetilde{B}$.
If $\widetilde{B}$ is an integer matrix whose principal part is
skew-symmetric, then it is 
easy to check that $\mu_k(\widetilde{B})$ is also an integer matrix 
with skew-symmetric principal part.

\subsection{}\label{modulemutation}
Let $T = T_1 \oplus \cdots \oplus T_r$ be a basic 
complete rigid $\la$-module
with $T_i$ indecomposable for all $i$.
Without loss of generality assume that $T_{r-n+1}, \ldots, T_r$ are
projective. 
Let $B(T) = (t_{ij})_{1 \le i,j \le r}$ be the $r \times r$-matrix defined
by
\[
t_{ij} = 
(\text{number of arrows $j \to i$ in $\Gamma_T$}) 
-(\text{number of arrows $i \to j$ in $\Gamma_T$}).
\]
Since the quiver $\Gamma_T$ does not have 2-cycles,
at least one of the two summands in the definition of $t_{ij}$
is zero.
Define 
$B(T)^\circ = (t_{ij})$
to be the $r \times (r-n)$-matrix obtained from
$B(T)$ by deleting the last $n$ columns.

For $k \in [1,r-n]$ there is a short exact sequence
\[
0 \to T_k \xrightarrow{f} \bigoplus_{t_{ik} > 0} T_i^{t_{ik}}
\to T_k^* \to 0
\]
where $f$ is a minimal left 
${\rm add}(T/T_k)$-approximation of $T_k$
(i.e. the map $\Hom_\la(f,T)$
is surjective, and every morphism $g$ with $gf = f$ is an 
isomorphism, see Section \ref{approxbasics} for a review of basic
results on approximations.)
Set 
\[
\mu_{T_k}(T) = T_k^* \oplus T/T_k.
\]
We show that $\mu_{T_k}(T)$ is again a basic complete rigid module
(Proposition \ref{compl1}).
In particular, $T_k^*$ is indecomposable.
We call $\mu_{T_k}(T)$ the 
{\it mutation of $T$ in direction $T_k$}.

Our second main result shows that the quivers of the endomorphism
algebras $\End_\la(T)$ and $\End_\la(\mu_{T_k}(T))$ are related
via Fomin and Zelevinsky's mutation rule:

\begin{Thm}\label{Mainresult1}
Let $\la$ be a preprojective algebra of Dynkin type
$\Delta$.
For a basic complete rigid $\la$-module $T$ as above and $k \in [1,r-n]$
we have
\[
B(\mu_{T_k}(T))^\circ = \mu_k(B(T)^\circ). 
\]
\end{Thm}

\subsection{}\label{deftdelta}
By $\Gamma_Q = (\Gamma_0,\Gamma_1,s,t,\tau)$ we denote the 
Auslander-Reiten quiver of $Q$, where $\tau$ is the Auslander-Reiten
translation.

Recall that
the vertices of $\Gamma_Q$ correspond to the isomorphism classes
of indecomposable $\field Q$-modules.
For an indecomposable $\field Q$-module $M$ let $[M]$ be its corresponding
vertex.
A vertex of $\Gamma_Q$ is called {\it projective} provided it corresponds
to an indecomposable projective $\field Q$-module.
For every indecomposable non-projective $\field Q$-module $X$ there
exists an Auslander-Reiten sequence
$$
0 \to \tau(X) \to \bigoplus_{i=1}^{n(X)} X_i  \to X \to 0.
$$
Here $\tau(X)$ and the $X_i$ are pairwise non-isomorphic indecomposable
$\field Q$-modules, which are uniquely dermined by $X$.
We also write $\tau([X]) = [\tau(X)]$.
The arrows of $\Gamma_Q$ are defined as follows:
Whenever 
there is an irreducible homomorphism $X \to Y$ between indecomposable
$\field Q$-modules $X$ and $Y$, then there is an arrow from $[X]$ to $[Y]$.
The Auslander-Reiten sequence above yields arrows $[\tau(X)] \to [X_i]$ and
$[X_i] \to [X]$, and all arrows of $\Gamma_Q$ are obtained in this way
from Auslander-Reiten sequences. 

In general, there can be more than one arrow between two 
vertices of an Auslander-Reiten quiver of an algebra.
But since we just work with path algebras of Dynkin quivers, this does
not occur in our situation.
For details on Auslander-Reiten sequences we refer to \cite{ARS} and
\cite{Ri1}.

It is well known (Gabriel's Theorem) that $\Gamma_Q$ has exactly $r$ 
vertices, say $1, \ldots, r$, and $n$ of these are projective.
Without loss of generality assume that the projective vertices are labelled 
by $r-n+1, \ldots, r$.
We define a new quiver $\Gamma_Q^*$ which is obtained from $\Gamma_Q$ 
by adding an arrow $x \to \tau(x)$ for each non-projective vertex $x$ 
of $\Gamma_Q$.
For the proof of Theorem \ref{introthm1}
we need the
following result \cite[Theorem 1]{GLS2}:

\begin{Thm}\label{startmoduletheorem}
There exists a basic complete rigid $\la$-module
$T_Q$ such that the quiver of the endomorphism algebra
$\End_\la(T_Q)$ is $\Gamma_Q^*$.
\end{Thm}

In comparison to \cite{GLS2} we changed our notation slightly: Here
we denote the
module ${\rm I}_{Q^{op}}$ which we constructed in \cite{GLS2} by $T_Q$.
Also we write $\Gamma_Q^*$ instead of $\check{A}_Q$.

Since $\la$ does not depend on the orientation of $Q$, we get that for every
Dynkin quiver $Q'$ of type $\Delta$ there is a basic complete rigid 
$\la$-module
$T_{Q'}$ such that the quiver of its endomorphism algebra
is $\Gamma_{Q'}^*$.

\subsection{}\label{clusterCN}
Note that in \cite{BFZ} only the cluster algebra structure
on $\C[G/N]$ is defined explicitly.
But one can easily modify it to get the one on $\C[N]$.
The cluster algebra structure on $\C[N]$ is defined
as follows, see \cite{GLS2} for details:
To any reduced word ${\ib}$ of the longest element $w_0$ of the
Weyl group, one can associate an initial seed
$(\tilde{\xb}',\widetilde{B}(\ib)')$ consisting of a list 
$\tilde{\xb}' = (\Delta(1,\ib)',\ldots,\Delta(r,\ib)')$ of $r$ 
distinguished elements of $\C[N]$ together with an
$r \times (r-n)$-matrix $\widetilde{B}({\ib})'$.
(Here we use the notation from \cite{GLS2}.)
This seed is described combinatorially in terms of $\ib$.
In particular, the elements of $\tilde{\xb}'$ are certain
explicit {\it generalized minors} attached to certain subwords of
${\ib}$. 
The other seeds are then produced inductively from
$(\tilde{\xb}',\widetilde{B}(\ib)')$ by a process of seed mutation
introduced by Fomin and Zelevinsky:

Assume 
$((f_1,\ldots,f_r),\widetilde{B})$ 
is a seed which was obtained  
by iterated seed mutation from our initial seed.
Thus the $f_i$ are certain elements in $\C[N]$ and 
$\widetilde{B} = (b_{ij})$ is a certain 
$r \times (r-n)$-matrix with integer entries.
For $1 \le k \le r-n$ define another $r \times (r-n)$-matrix 
$\mu_k(\widetilde{B}) = (b_{ij}')$ as in Section \ref{matmutation} above, 
and let
$$
f_k' = 
\frac{\prod_{b_{ik}>0} f_i^{b_{ik}} + \prod_{b_{ik}<0} f_i^{-b_{ik}}}{f_k}.
$$
The choice of the initial seed ensures that $f_k'$ is again an element
in $\C[N]$.
This follows from the results in \cite{BFZ}.
Let $\mu_k(f_1,\ldots,f_r)$ be the $r$-tuple obtained from $(f_1,\ldots,f_r)$
by replacing the entry $f_k$ by $f_k'$.
Then 
$
(\mu_k(f_1,\ldots,f_r),\mu_k(\widetilde{B}))
$ 
is the seed obtained from
$((f_1,\ldots,f_r),\widetilde{B})$ by mutation in direction $k$.
Thus, starting with our initial seed, an inductive 
combinatorial procedure gives all the other seeds.
Each seed has $r-n$ neighbouring seeds obtained by just one mutation.

Suppose that ${\ib}$ is {\it adapted} to the quiver $Q$, 
in the sense of \cite{Lu}. 
Then it is shown in 
\cite{GLS2} that the matrix $\widetilde{B}({\ib})'$ coincides with the matrix
$B(T_Q)^\circ$,
where $T_Q$ is the basic complete rigid module of 
Theorem~\ref{startmoduletheorem}.
Moreover, if we write $T_Q=T_1\oplus \cdots \oplus T_r$
and $x_i = \varphi_{T_i} \ (1\le i\le r)$ (see \ref{lift} above), then 
$\tilde{\xb}'$ coincides with $(x_1, \ldots , x_r)$.
In other words, the rigid module $T_Q$ can be regarded as 
a lift to $\md(\Lam)$ of the initial seed $(\tilde{\xb}',\widetilde{B}(\ib)')$.

Let $\T_\Lam$ be the graph with vertices the isomorphism
classes of basic maximal rigid $\Lam$-modules and with edges given 
by mutations. Let $\T_\Lam^\circ$ denote the connected 
component of $\T_\Lam$ containing $T_Q$.
To a vertex $R=R_1\oplus \cdots \oplus R_r$ we attach the 
$r$-tuple of regular functions 
$\xb(R) = (\varphi_{R_1},\ldots,\varphi_{R_r})$.

For a rigid $\Lam$-module $M$ 
the closure of the orbit $\orb_M$ in the corresponding
module variety is an irreducible component, see Section \ref{preproj} 
(in particular Corollary \ref{corlu2}).
Therefore $M$ is a generic point of the irreducible component 
$\overline{\orb_M}$, see Section \ref{semcan}.
This implies that $\xb(R)$ is a collection of elements of the dual
semicanonical basis $\SB^*$. 
Define $B(R)^\circ$ as in Section \ref{modulemutation}.

Let $\G$ be the {\it exchange graph} of the cluster algebra
$\C[N]$, that is, the graph with vertices the seeds and
edges given by seed mutation.

Using Theorem~\ref{Mainresult1} and the multiplication formula
for functions $\varphi_M$ of \cite{GLS3}, we can then deduce
the following theorem.

\begin{Thm}\label{clustermonom}
With the notation above the following hold:
\begin{itemize}

\item[(1)]
For each vertex $R$ of $\T_\Lam^\circ$, 
$\xb(R)$ is a cluster of $\C[N]$;

\item[(2)]
The map $R \mapsto (\xb(R),B(R)^\circ)$ induces an isomorphism of graphs
from $\T_\Lam^\circ$ to $\G$;

\item[(3)]
The cluster monomials belong to 
the dual semicanonical basis $\SB^*$. 
\end{itemize}
\end{Thm}

Note that this proves in particular that the cluster monomials 
are linearly independent, which is not obvious from 
their definition.
On the other hand, with the exception of Dynkin types
$\A_n$ with $n \le 4$, it is known that there exist irreducible
components of varieties of $\Lam$-modules without rigid modules.
Therefore, it also follows from Theorem~\ref{clustermonom} that
when $\C[N]$ is a cluster algebra of infinite
type, the cluster monomials form a proper subset of $\SB^*$
and do not span $\C[N]$.

It is also worth noting that there is no known 
algorithm to calculate the semicanonical basis
(in \cite{Lu1}, a similar basis of constructible
functions for the group algebra of a Weyl group
is considered to be ``probably uncomputable''). 
Therefore it is remarkable that, by Theorem~\ref{clustermonom}, 
a large family of elements of $\SB^*$ can be 
evaluated by a combinatorial algorithm, namely
by repeated applications of the exchange formula
for cluster mutation.

We conjecture that the graph $\T_\Lam$ is connected, so that 
Theorem~\ref{clustermonom} should hold with $\T_\Lam$ instead
of $\T_\Lam^\circ$.

\subsection{}
The paper is organized as follows:

In Section \ref{preliminary} we recall mostly known results
from the representation theory of algebras.
Section \ref{maxorthogonal} discusses maximal 1-orthogonal
modules, which were recently introduced and studied by Iyama
\cite{I}, \cite{I1}.
We repeat some of Iyama's results, and then study the functor
$\Hom_A(-,T)$ associated to a maximal 1-orthogonal module $T$.
In Section \ref{mutationsection} we introduce a mutation operation on 
basic rigid $\la$-modules and study the corresponding exchange
sequences.

The endomorphism algebras $\End_\la(T)$ of maximal
rigid $\la$-modules are studied in detail in Section \ref{endosection},
which contains a proof of Theorem \ref{introthm1}.

Next, in Section \ref{graphmutation} we use the results from 
Section \ref{endosection} to prove Theorem \ref{Mainresult1}.
Some examples are given in Section \ref{example}.

In Section \ref{cluster} we recall some backround results
on cluster algebras and semicanonical bases, and
we prove Theorem~\ref{clustermonom}.


\section{Preliminary results}\label{preliminary}


In this section let $A$ be a finite-dimensional algebra over an 
algebraically closed field $\field$.
For backround material concerning the representation theory of
finite-dimensional algebras we refer to \cite{ARS} and
\cite{Ri1}.

\subsection{Approximations of modules}\label{approxbasics}
We recall some well known results from the theory of approximations
of modules.
Let $M$ be an $A$-module.
A homomorphism $f\colon X \to M'$ in $\md(A)$ is 
a {\it left ${\rm add}(M)$-approximation} of $X$ if $M' \in {\rm add}(M)$ 
and the induced map
\[
\Hom_A(f,M)\colon \Hom_A(M',M) \to \Hom_A(X,M) 
\]
(which maps a homomorphism $g\colon M' \to M$ to $gf\colon X \to M$)
is surjective,
i.e. every homomorphism $X \to M$ 
factors through $f$.

A morphism $f\colon V \to W$ is called
{\it left minimal} if every morphism $g\colon W \to W$ with $gf = f$ is an 
isomorphism.

Dually, one defines right ${\rm add}(M)$-approximations and
right minimal morphisms:
A homomorphism $f\colon M' \to X$ is a {\it right ${\rm add}(T)$-approximation}
of $X$ if $M' \in {\rm add}(M)$ 
and the induced map
\[
\Hom_A(M,f)\colon \Hom_A(M,M') \to \Hom_A(M,X) 
\]
(which maps a homomorphism $g\colon M \to M'$ to $fg\colon M \to X$)
is surjective,
i.e. every homomorphism $M \to X$ 
factors through $f$.
A morphism $f\colon V \to W$ is called
{\it right minimal} if every morphism $g\colon V \to V$ with $fg = f$ is an 
isomorphism.

The following results are well known:

\begin{Lem}\label{approx1}
Let $f\colon X \to M$ be a homomorphism of $A$-modules.
Then there exists a decomposition $M = M_1 \oplus M_2$
with $\Ima(f) \subseteq M_1$ and $f_1\colon X \to M_1$ 
is left minimal, where $f_1 = \pi_1f$ with
$\pi_1\colon M \to M_1$ the canonical projection.
\end{Lem}

\begin{Lem}\label{approx2}
Let 
\[
0 \to X \xrightarrow{f} Y \xrightarrow{g} 
Z \to 0
\]
be a non-split short exact sequence of $A$-modules.
Then the following hold:
\begin{itemize}

\item
If $Z$ is indecomposable, then $f$ is left minimal;

\item
If $X$ is indecomposable, then $g$ is right minimal.

\end{itemize}
\end{Lem}

\begin{Lem}\label{approx3}
If $f_i\colon X \to M_i$ $(i=1,2)$ are minimal left 
${\rm add}(M)$-approximations,
then $M_1 \cong M_2$ and $\Coker(f_1) \cong \Coker(f_2)$. 
\end{Lem}

\begin{Lem}\label{approx5}
If $f_i\colon X \to M_i$ $(i=1,2)$ are left ${\rm add}(M)$-approximations
with $M_1 \cong M_2$, then $\Coker(f_1) \cong \Coker(f_2)$. 
\end{Lem}

There are obvious duals of Lemmas \ref{approx1}, \ref{approx3} and
\ref{approx5}.

\begin{Lem}\label{approx4}
Let
\[
0 \to X \xrightarrow{f} T' \xrightarrow{g} Y \to 0
\]
be a non-split short exact sequence 
with $T' \in {\rm add}(T)$ for some module $T$.
Then the following hold:
\begin{itemize}

\item
If $\Ext_A^1(Y,T) = 0$, then $f$ is a left ${\rm add}(T)$-approximation;

\item
If $\Ext_A^1(T,X) = 0$, then $g$ is a right ${\rm add}(T)$-approximation.

\end{itemize}
\end{Lem}

\begin{proof}
Applying $\Hom_A(-,T)$ we get an exact sequence
$$
0 \to \Hom_A(Y,T) \to \Hom_A(T',T) \xrightarrow{\Hom_A(f,T)} \Hom_A(X,T)
\to \Ext_A^1(Y,T) = 0.
$$
Thus $\Hom_A(f,T)$ is surjective, i.e. $f$ is a left 
${\rm add}(T)$-approximation.
If we apply the functor
$\Hom_A(T,-)$ we obtain the corresponding dual result for $g$.
\end{proof}

\subsection{Endomorphism algebras and their quivers}\label{endoquiver}
Let $M = M_1^{n_1} \oplus \cdots \oplus M_t^{n_t}$ 
be a finite-dimensional $A$-module, where the $M_i$
are pairwise non-isomorphic indecomposable
modules and $n_i \ge 1$.
As before let $S_i = S_{M_i}$ be the simple $\End_A(M)$-module
corresponding to $M_i$, and let
$P_i = \Hom_A(M_i,M)$ be the indecomposable projective $\End_A(M)$-module with
top $S_i$.
For $1 \le i,j \le t$ the following numbers are equal:
\begin{itemize}

\item
The number of arrows $i \to j$ in the quiver $\Gamma_M$ of
$\End_A(M)$;

\item
$\dm \Ext_{\End_A(M)}^1(S_i,S_j)$;

\item
The dimension of the space 
of irreducible maps $M_i \to M_j$ in the category ${\rm add}(M)$;

\item
The dimension of the space 
of irreducible maps $P_j \to P_i$ in the category  
${\rm add}(P_1 \oplus \cdots \oplus P_t)$
of projective $\End_A(M)$-modules.

\end{itemize}
Furthermore, let $f\colon M_i \to M'$ (resp. $g\colon M'' \to M_i$) 
be a minimal left (resp. right) ${\rm add}(M/(M_i^{n_i}))$-approximation
of $M_i$.
If $i \not= j$, then we have
\begin{align*}
\dm \Ext_{\End_A(M)}^1(S_i,S_j) &= [M':M_j],\\
\dm \Ext_{\End_A(M)}^1(S_j,S_i) &= [M'':M_j].
\end{align*}
The above facts are important and will be often used, in particular 
in Sections
\ref{endosection} and
\ref{graphmutation}.

\subsection{Tilting modules and derived equivalences}\label{tiltingsection}
An $A$-module $T$ 
is a {\it tilting module} if the following three conditions hold:
\begin{itemize}

\item[(1)] $\pdim(T) \le 1$;

\item[(2)]
$\Ext_A^1(T,T) = 0$;

\item[(3)] There exists a short exact sequence
\[
0 \to A \to T' \to T'' \to 0
\]
with $T',T'' \in {\rm add}(T)$.

\end{itemize}

By $D^b(A)$ we denote the derived category of bounded complexes
of $A$-modules.
Let $B$ be another finite-dimensional $\field$-algebra.
The algebras $A$ and $B$ are {\it derived equivalent} if the categories
$D^b(A)$ and $D^b(B)$ are equivalent as triangulated categories, see
for example \cite[Section 0]{H}.

\begin{Thm}[{\cite[Section 1.7]{H}}]\label{Happel}
If $T$ is a tilting module over $A$, then
$A$ and $\End_A(T)^\op$ are derived equivalent.
\end{Thm}

The following theorem is also well known:

\begin{Thm}[{\cite[Section 1.4]{H}}]
If $A$ and $B$ are derived equivalent,
then $\gldim(A) < \infty$ if and only if $\gldim(B) < \infty$.
\end{Thm}

For details on tilting theory and derived categories we refer to
\cite{HR}, \cite{H}.

\subsection{Loops and 2-cycles}
For the rest of this section we assume 
that $A = \field \Gamma/I$, where $\Gamma$ is a finite quiver and
$I$ is an admissible ideal.
We denote the simple $A$-module corresponding to a vertex $l$ of $\Gamma$
by $S_l$.
The path of length 0 at $l$ is denoted by $e_l$.
Thus $e_p \field \Gamma e_q$ is the vector space with
basis the set of paths in $\Gamma$ which start at $q$ and end at $p$.

The following result is proved in \cite{Ig}:

\begin{Thm}[Igusa]\label{Igusa}
If $\gldim(A) < \infty$,
then the quiver of $A$ has no loops.
\end{Thm}

Let $[A,A]$ be the {\it commutator subgroup} of $A$.
This is the subspace of $A$ generated by all commutators
$[a,b] = ab-ba$ with $a,b \in A$.
The following result can be found in  \cite[Satz 5]{Le}:

\begin{Thm}[Lenzing]\label{nilpotent}
If $\gldim(A) < \infty$, then every nilpotent element 
$a \in A$ lies in $[A,A]$.
\end{Thm}

As mentioned in \cite{Ig}, 
Theorem \ref{Igusa} is already contained implicitly
in Theorem \ref{nilpotent}.
The next lemma follows directly from \cite{Bo}, see also \cite{BK}:

\begin{Lem}\label{2cyclelemma}
Let $w$ be a path of length two in $e_p \field \Gamma e_q$.
If $w - c \in I$ for some $c \in e_p \field \Gamma e_q$ with $w \not= c$, 
then $\Ext_A^2(S_q,S_p) \not= 0$.
\end{Lem}

From Theorem \ref{nilpotent} we can deduce the following result:

\begin{Prop}\label{2cycleresult}
Assume that $\gldim(A) < \infty$ and that
the quiver of $A$ has a 2-cycle. 
Then $\Ext_A^2(S,S) \not= 0$ 
for some simple $A$-module $S$.
\end{Prop}

\begin{proof}
For some vertices $i \not= j$ in $\Gamma$ there are
arrows $a\colon j \to i$ and $b\colon i \to j$ in $\Gamma$.
We have $ab \in e_i \field \Gamma e_i$.
We claim that 
$\Ext_A^2(S_i,S_i) \not= 0$ or
$\Ext_A^2(S_j,S_j) \not= 0$:

We know that $ab$ is nilpotent in $A$.
Thus by Theorem \ref{nilpotent} we have $ab \in [A,A]$.
One easily checks that every commutator is a linear combination
of commutators of the form $[v,w]$ where $v$ and $w$ are paths
in $\Gamma$.
Thus 
\[
ab - \sum_{l=1}^m \lambda_l [v_l,w_l] \in I
\]
for some appropriate paths $v_l$ and $w_l$ in $\Gamma$ and
some scalars $\lambda_l \in \field$.
Without loss of generality we can assume the following:
\begin{itemize}

\item
$v_1 = a$ and $w_1 = b$;

\item
$\{ v_s,w_s \} \not= \{ v_t,w_t \}$
for all $s \not= t$;

\item
$\{ v_lw_l, w_lv_l \} \cap \{ ab,ba \} = \emptyset$ for
all $2 \le l \le m$. 

\end{itemize}
Note that we do not assume that the scalars $\lambda_l$ are all non-zero.
Set 
\[
C = ab - \sum_{l=1}^m \lambda_l [v_l,w_l].
\]
We get 
\[
e_iCe_i = ab - \lambda_1 ab - \sum_{l=2}^m \lambda_l e_i[v_l,w_l]e_i \in I.
\]
If $\lambda_1 \not= 1$, then Lemma \ref{2cyclelemma} applies and we get
$\Ext_A^2(S_i,S_i) \not= 0$.
If $\lambda_1 \not= 0$, then
\[
e_jCe_j = \lambda_1 ba - \sum_{l=2}^m \lambda_l e_j[v_l,w_l]e_j \in I.
\]
Again we can apply Lemma \ref{2cyclelemma} to this situation and
get $\Ext_A^2(S_j,S_j) \not= 0$.
This finishes the proof.
\end{proof}

\subsection{Preprojective algebras}\label{preproj}
We recall some results on preprojective algebras.
There is a symmetric bilinear form $(-,-): \Z^n \times \Z^n \to \Z$
associated to $Q$
defined by
\[
(d,e) = 2 \sum_{i \in Q_0} d_ie_i -
\sum_{\alpha \in Q_1} (d_{s(\alpha)}e_{t(\alpha)} + 
e_{s(\alpha)}d_{t(\alpha)}).
\]
The dimension vector of a module $M$ is denoted by $\dimv(M)$.
The following lemma is due to Crawley-Boevey \cite[Lemma 1]{CB}.

\begin{Lem}\label{extandhom}
For $\la$-modules $X$ and $Y$ we have
\[
\dm \Ext_\la^1(X,Y) = \dm \Hom_\la(X,Y) + \dm \Hom_\la(Y,X) - 
(\dimv(X),\dimv(Y)).
\]
\end{Lem}

\begin{Cor}
$\dm \Ext_\la^1(X,X)$ is even, and 
$\dm \Ext_\la^1(X,Y) = \dm \Ext_\la^1(Y,X)$.
\end{Cor}

Let $\beta = (\beta_1, \ldots, \beta_n) \in \N^n$.
By 
\[
\rep(Q,\beta) = \prod_{\alpha \in Q_1} 
\Hom_{\field}(\field^{\beta_{s(\alpha)}},
\field^{\beta_{t(\alpha)}}) 
\]
we denote the affine space of representations of $Q$ with dimension
vector $\beta$.
Furthermore, let $\la_\beta$ be the affine variety of elements
\[
(f_\alpha,f_{\alpha^*})_{\alpha \in Q_1} \in
\prod_{\alpha \in Q_1} \left( \Hom_{\field}(\field^{\beta_{s(\alpha)}},
\field^{\beta_{t(\alpha)}}) \times \Hom_{\field}(\field^{\beta_{t(\alpha)}},
\field^{\beta_{s(\alpha)}}) \right) 
\]
such that for all $i \in Q_0$ we have
\[ 
\sum_{\alpha \in Q_1: s(\alpha) = i} f_{\alpha^*}f_\alpha =
\sum_{\alpha \in Q_1: t(\alpha) = i} f_\alpha f_{\alpha^*}.
\]
Note that $\rep(Q,\beta)$ can be considered as a subvariety of $\la_\beta$.
Since $Q$ is a Dynkin quiver, the variety $\la_\beta$ coincides with
Lusztig's nilpotent variety associated to $Q$.
One can interpret $\la_\beta$ as the variety of $\la$-modules with dimension
vector $\beta$.
The group 
\[
\GL_\beta = \prod_{i=1}^n {\rm GL}_{\beta_i}(\field)
\]  
acts by conjugation
on $\rep(Q,\beta)$ and $\la_\beta$:

For 
$g = (g_1,\ldots,g_n) \in \GL_\beta$ and
$x = (f_\alpha,f_{\alpha^*})_{\alpha \in Q_1} \in \la_\beta$ define
\[
g \cdot x = (g_{t(\alpha)}f_\alpha g_{s(\alpha)}^{-1},
g_{s(\alpha)}f_{\alpha^*} g_{t(\alpha)}^{-1})_{\alpha \in Q_1}.
\]
The action on $\rep(Q,\beta)$ is obtained via restricting the action on
$\la_\beta$ to $\rep(Q,\beta)$.

The isomorphism classes of representations in $\rep(Q,\beta)$ and $\la$-modules
in $\la_\beta$, respectively, correspond to the orbits of these actions.
The $\GL_\beta$-orbit of some $M \in \rep(Q,\beta)$ or $M \in \la_\beta$ is 
denoted by $\orb_M$.
Since $Q$ is a Dynkin quiver, there are only finitely many
orbits in $\rep(Q,\beta)$.
Let 
$
\pi_\beta\colon \la_\beta \to \rep(Q,\beta)
$
be the canonical projection morphism.
Lusztig shows that
$
\orb \mapsto \overline{\pi_\beta^{-1}(\orb)}
$
defines a one-to-one correspondence between the $\GL_\beta$-orbits
in $\rep(Q,\beta)$ and the set $\irr(\la_\beta)$ of irreducible components
of $\la_\beta$.
He also proved that 
all irreducible components of $\la_\beta$ have dimension
\[
\sum_{\alpha \in Q_1} \beta_{s(\alpha)}\beta_{t(\alpha)},
\]
see \cite[Section 12]{Lu}.
For a $\GL_\beta$-orbit $\orb$ in $\la_\beta$ let 
${\rm codim}\, \orb = \dm \la_\beta - \dm \orb$ be its codimension.

\begin{Lem}\label{corlu1}
For any $\la$-module $M$ we have
$
\dm \Ext_\la^1(M,M) = 2\, {\rm codim}\, \orb_M$.
\end{Lem}

\begin{proof}
Set $\beta = \dimv(M)$.
By Lemma \ref{extandhom} we have
$
\dm \Ext_\la^1(M,M) = 2\, \dm \End_\la(M) - (\beta,\beta)$.
Furthermore, 
$
\dm \orb_M = \dm \GL_\beta - \dm \End_\la(M)$.
Thus
\[
{\rm codim}\, \orb_M = \dm \la_\beta - \dm \orb_M =
\sum_{\alpha \in Q_1} \beta_{s(\alpha)}\beta_{t(\alpha)} - 
\sum_{i=1}^n \beta_i^2
+ \dm \End_\la(M).
\]
Combining these equations yields the result.
\end{proof}

\begin{Cor}\label{corlu2}
For a $\la$-module $M$ with dimension vector
$\beta$ the following are equivalent:
\begin{itemize}

\item
The closure $\overline{\orb_M}$ of $\orb_M$ is an irreducible
component of $\la_\beta$;

\item
The orbit $\orb_M$ is open in $\la_\beta$;

\item
$\Ext_\la^1(M,M) = 0$.

\end{itemize}
\end{Cor}


\section{Maximal $1$-orthogonal modules}
\label{maxorthogonal}

In this section let $A$ be a finite-dimensional algebra over an 
algebraically closed field $\field$.

\subsection{Iyama's results}\label{iyamaresults}
We recall some of Iyama's recent results on maximal 1-orthogonal
modules \cite{I}, \cite{I1}.

\begin{Thm}[{\cite[Theorem 5.1 (3)]{I1}}]\label{Iyama1}
Let $T$ be a rigid $A$-module and assume that
$T$ is a generator-cogenerator.
If $\gldim(\End_A(T)) \le 3$, then $T$ is maximal $1$-orthogonal.
\end{Thm}

\begin{Thm}[{\cite[Theorem 0.2]{I1}}]\label{Iyama2}
If $T$ is a maximal $1$-orthogonal $A$-module, then
\[
\gldim(\End_A(T)) \le 3 \text{ and } \domdim(\End_A(T)) \ge 3.
\]
\end{Thm}

\begin{Thm}[{\cite[Theorem 5.3.2]{I1}}]\label{Iyama3}
Let $T_1$ and $T_2$ be maximal $1$-orthogonal $A$-modules.
Then
$T = \Hom_A(T_2,T_1)$ is a tilting module over
$\End_A(T_1)$, and we have
\[
\End_{\End_A(T_1)}(T) \cong \End_A(T_2)^\op.
\]
In particular, $\End_A(T_1)$ and $\End_A(T_2)$ are derived equivalent,
and $\Sigma(T_1) = \Sigma(T_2)$.
\end{Thm}

\subsection{A Hom-functor}

\begin{Prop}\label{homfunctor}
Let $T$ be a maximal 1-orthogonal $A$-module.
Then the contravariant functor 
\[
F_T = \Hom_A(-,T)\colon \md(A) \to \md(\End_A(T))
\]
yields an anti-equivalence of categories
\[
\md(A) \to {\mathcal P}(\End_A(T))
\]
where ${\mathcal P}(\End_A(T)) \subset \md(\End_A(T))$ 
denotes the full subcategory of all
$\End_A(T)$-modules of projective dimension at most one.
\end{Prop}

\begin{proof}
Set $E = \End_A(T)$.
Let $X \in \md(A)$, and let $f\colon X \to T'$ be a minimal left
${\rm add}(T)$-approximation of $X$.
Since $T$ is a cogenerator of $\md(A)$, we know that $f$ is a
monomorphism.
Thus there is a short exact sequence
\[
0 \to X \xrightarrow{f} T' \to T'' \to 0.
\]  
Applying $\Hom_A(-,T)$ to this sequence
yields an exact sequence
\begin{multline*}
0 \to \Hom_A(T'',T) \to \Hom_A(T',T) \xrightarrow{\Hom_A(f,T)} 
\Hom_A(X,T)\\ 
\to \Ext_A^1(T'',T)
\to \Ext_A^1(T',T).
\end{multline*}
Since $T$ is rigid and $T' \in {\rm add}(T)$, we get $\Ext_A^1(T',T) = 0$.
The map $f$ is an ${\rm add}(T)$-approximation, thus $\Hom_A(f,T)$ is
surjective.
This implies $\Ext_A^1(T'',T) = 0$.
But $T$ is maximal 1-orthogonal, which yields $T'' \in {\rm add}(T)$.
Thus the exact sequence
\[
0 \to \Hom_A(T'',T) \to \Hom_A(T',T) \xrightarrow{\Hom_A(f,T)} 
\Hom_A(X,T) \to 0
\]
is a projective resolution of $F_T(X) = \Hom_A(X,T)$, and we conclude
$$
\pdim_E(F_T(X)) \le 1.
$$

Now the dual of \cite[Lemma 1.3(b)]{APR} yields that $F_T$
induces a fully faithful functor $\md(A) \to {\mathcal P}(E)$.

Finally, we show that every object $Y \in {\mathcal P}(E)$ lies in the 
image of $F_T$.
There exists a projective resolution of the form
\[
0 \to \Hom_A(T'',T) \xrightarrow{G} \Hom_A(T',T) \to Y \to 0
\]
with $T',T'' \in {\rm add}(T)$.
Since $\Hom_A(T',T)$ and $\Hom_A(T'',T)$ are projective,
there exists some $g \in \Hom_A(T',T'')$ such that $G = \Hom_A(g,T)$.
We claim that $g$ must be surjective.
Let us assume that $g$ is not surjective.
Then there exists a non-zero homomorphism $h \in \Hom_A(T'',T)$ such that
$hg = 0$.
Here we use that $T$ is a cogenerator of $\md(A)$.
This implies $G(h) = 0$, which is a contradiction to $G$ being injective.

Thus there is a short exact sequence
\[
0 \to \Ker(g) \to T' \xrightarrow{g} T'' \to 0.
\]
Applying $\Hom_A(-,T)$ to this sequence yields an exact sequence
\[
0 \to \Hom_A(T'',T) \xrightarrow{G} \Hom_A(T',T) \to \Hom_A(\Ker(g),T) \to 0.
\]
Here we used that $\Ext_A^1(T'',T) = 0$.
Thus we get $Y \cong \Hom_A(\Ker(g),T)$.
This finishes the proof of Proposition \ref{homfunctor}.
\end{proof}

A more general result than the above Proposition \ref{homfunctor}
can be found in \cite[Theorem 5.3.4]{I1}.

Let $F_T$ be as in Proposition \ref{homfunctor}.
Then $F_T$ has the following properties:
\begin{itemize}

\item
If $X$ is an indecomposable $A$-module, then $F_T(X)$ is indecomposable;

\item
$F_T$ reflects isomorphism classes, i.e.
if $F_T(X) \cong F_T(Y)$ for some $A$-modules $X$ and $Y$, then
$X \cong Y$.

\end{itemize}
Note that the functor $F_T$ is not exact.

The following proposition is due to Iyama:

\begin{Prop}\label{numberofcomplements1}
Let $T$ be a basic maximal 1-orthogonal $A$-module, and let
$X$ be an indecomposable direct summand of $T$.
Then there is at most one indecomposable $A$-module 
$Y$ such that $X \not\cong Y$ and $Y \oplus T/X$ is maximal 1-orthogonal.
\end{Prop}

\begin{proof}
Assume there are two non-isomorphic indecomposable $A$-modules
$Y_i$ where $i = 1,2$ such that $X \not\cong Y_i$ and 
$Y_i \oplus T/X$ is maximal 1-orthogonal. 
By Theorem \ref{Iyama3} we know that $\Hom_A(Y_i \oplus T/X,T)$
is a tilting module over $\End_A(T)$.

Thus the almost complete tilting module $\Hom_A(T/X,T)$ 
over $\End_A(T)$ has three indecomposable complements,
namely $\Hom_A(X,T)$, $\Hom_A(Y_1,T)$ and $\Hom_A(Y_2,T)$.
Here we use that the functor $F_T$ as defined in
Proposition \ref{homfunctor} preserves indecomposables and reflects
isomorphism classes.
By \cite[Proposition 1.3]{RS2} an almost
complete tilting module has at most two indecomposable complements,
a contradiction.
\end{proof}


\section{Mutations of rigid modules}\label{mutationsection}


\subsection{}
In this section we work with modules over the 
preprojective algebra $\la$.

\begin{Lem}\label{prop1}
Let $T$ and $X$ be rigid $\la$-modules.
If
\[
0 \to X \xrightarrow{f} T' \xrightarrow{g} Y \to 0
\]
is a short exact sequence with $f$ a left ${\rm add}(T)$-approximation, then
$T \oplus Y$ is rigid.
\end{Lem}

\begin{proof}
First, we prove that $\Ext_\la^1(Y,T) = 0$:
We apply $\Hom_\la(-,T)$ and get an exact sequence
\begin{multline*}
0 \to \Hom_\la(Y,T) \to \Hom_\la(T',T) 
\xrightarrow{\Hom_\la(f,T)} \Hom_\la(X,T)\\ 
\to \Ext_\la^1(Y,T) 
\to \Ext_\la^1(T',T) = 0.
\end{multline*}
Since $f$ is a left ${\rm add}(T)$-approximation, we know that $\Hom_\la(f,T)$
is surjective.
Thus $\Ext_\la^1(Y,T) = 0$.
Next, we show that $\Ext_\la^1(Y,Y) = 0$.
This is similar to the proof of \cite[Lemma 6.7]{BMRRT}.
We apply $\Hom_\la(X,-)$ and get an exact sequence
$$
0 \to \Hom_\la(X,X) \to \Hom_\la(X,T') \xrightarrow{\Hom_\la(X,g)} 
\Hom_\la(X,Y) \to \Ext_\la^1(X,X) = 0.
$$
Thus $\Hom_\la(X,g)$ is surjective, i.e. every morphism $h\colon X \to Y$
factors through $g\colon T' \to Y$. 
We have $\Ext_\la^1(Y,T) = 0$, and by 
Lemma \ref{extandhom} we get $\Ext_\la^1(T,Y) = 0$.

Applying $\Hom_\la(-,Y)$ yields an exact sequence
\begin{multline*}
0 \to \Hom_\la(Y,Y) \to \Hom_\la(T',Y) \xrightarrow{\Hom_\la(f,Y)} 
\Hom_\la(X,Y)\\ 
\to \Ext_\la^1(Y,Y) 
\to \Ext_\la^1(T',Y) = 0.
\end{multline*}
To show that $\Ext_\la^1(Y,Y) = 0$ it is enough to show that $\Hom_\la(f,Y)$
is surjective, i.e. that every map $h\colon X \to Y$ factors through 
$f\colon X \to T'$. 
Since $\Hom_\la(X,g)$ is surjective, there is a morphism $t\colon X \to T'$
such that $gt = h$.
Since $f$ is a left ${\rm add}(T)$-approximation, there is a morphism
$s\colon T' \to T'$ such that $sf = t$.
So
\[
h = gt = gsf. 
\]
Thus $h$ factors through $f$, which implies $\Ext_\la^1(Y,Y) = 0$.
\end{proof}

\begin{Cor}\label{cor1}
Let $T$ and $X$ be rigid $\la$-modules.
If $T$ is maximal rigid, then there exists a short exact sequence
\[
0 \to X \to T' \to T'' \to 0
\]
with $T',T'' \in {\rm add}(T)$.
\end{Cor}

\begin{proof}
In the situation of Lemma \ref{prop1}, the maximality of
$T$ implies $Y \in {\rm add}(T)$.
\end{proof}

Note that there exist dual results for Lemma \ref{prop1} and
Corollary \ref{cor1}, involving right instead of left 
${\rm add}(T)$-approximations.

\begin{Cor}\label{projdim1}
Let $T$ and $X$ be rigid $\la$-modules, and set $E = \End_\la(T)$.
If $T$ is maximal rigid, then 
\[
\pdim_E(\Hom_\la(X,T)) \le 1.
\]
\end{Cor}

\begin{proof}
Apply $\Hom_\la(-,T)$ to the short exact sequence in Corollary \ref{cor1}.
This yields a projective resolution
\[
0 \to \Hom_\la(T'',T) \to \Hom_\la(T',T) \to \Hom_\la(X,T) \to 0.
\]
\end{proof}

\begin{Thm}\label{tilting}
For $i=1,2$ let $T_i$ be a maximal rigid $\la$-module,
let $E_i = \End_\la(T_i)$, 
and set 
$T = \Hom_\la(T_2,T_1)$.
Then $T$ is a tilting module over $E_1$ and
we have
\[
\End_{E_1}(T) \cong E_2^\op.
\]
In particular, $E_1$ and $E_2$ are derived equivalent,
and $\Sigma(T_1) = \Sigma(T_2)$.
\end{Thm}

\begin{proof}
It follows from Corollary \ref{cor1} 
that there exists a short exact sequence
\[
0 \to T_2 \to T_1' \to T_1'' \to 0
\]
with $T_1',T_1'' \in {\rm add}(T_1)$.
Then Corollary \ref{projdim1} yields that
the projective dimension of $T = \Hom_\la(T_2,T_1)$ regarded as 
an $E_1$-module is at most one.
This shows that $T$ satisfies condition (1) in the definition of a 
tilting module.

Now the rest of the proof is identical to Iyama's proof of Theorem 
\ref{Iyama3}.
\end{proof}

The above theorem is very similar to Iyama's Theorem \ref{Iyama3}.
However there are two differences:
Theorem \ref{Iyama3} holds for arbitrary finite-dimensional algebras, 
whereas we restrict to the class of finite-dimensional
preprojective algebras.
On the other hand, the modules $T_i$ in Theorem \ref{Iyama3} are
assumed to be maximal 1-orthogonal, which is stronger than our assumption
that the $T_i$ are maximal rigid.

In both cases, one has to prove that $\pdim(T) \le 1$.
This is the only place in the proof of Theorem \ref{Iyama3}
where the maximal 1-orthogonality of the $T_i$ is needed.
For the rest of his proof Iyama just needs that the $T_i$ are
maximal rigid.
If we restrict now to finite-dimensional preprojective algebras, then one
can prove that $\pdim(T) \le 1$ under the weaker assumption that the $T_i$
are maximal rigid.

Later on we then use Theorem \ref{tilting} in order to show 
that
the maximal rigid $\la$-modules are in fact maximal 1-orthogonal 
and also complete rigid.

\begin{Cor}\label{maxiscomplete}
For a $\la$-module $M$ the following are equivalent:
\begin{itemize}

\item
$M$ is maximal rigid;

\item
$M$ is complete rigid.

\end{itemize}
\end{Cor}

\begin{proof}
By Theorem \ref{upperbound} every complete rigid module is
maximal rigid.
For the other direction, let
$T_1$ be a complete rigid $\la$-module.
Such a module exists by Theorem \ref{startmoduletheorem}.
Let $T_2$ be a maximal rigid module.
Theorem \ref{tilting} implies that
$\End_\la(T_1)$ and $\End_\la(T_2)$ are derived equivalent.
Thus $\Sigma(T_1) = \Sigma(T_2) = r$, in other words, $T_2$ is also complete
rigid.
\end{proof}

\begin{Prop}\label{compl1}
Let $T \oplus X$ be a basic rigid $\la$-module such that the following
hold:
\begin{itemize}

\item $X$ is indecomposable;

\item $\la \in {\rm add}(T)$.

\end{itemize}
Then there exists a short exact sequence
\[
0 \to X \xrightarrow{f} T' \xrightarrow{g} Y \to 0
\]
such that the following hold:
\begin{itemize}

\item
$f$ is a minimal left ${\rm add}(T)$-approximation;

\item
$g$ is a minimal right ${\rm add}(T)$-approximation;

\item $T \oplus Y$ is basic rigid;

\item $Y$ is indecomposable and $X \not\cong Y$.

\end{itemize}
\end{Prop}

\begin{proof}
Let $f\colon X \to T'$ be a minimal left ${\rm add}(T)$-approximation of
$X$.
Since $\la \in {\rm add}(T)$, we know that $f$ is a monomorphism.
Let $Y$ be the cokernel of $f$ and let $g\colon T' \to Y$ be the 
projection map.
Thus we have a short exact sequence
\[
0 \to X \xrightarrow{f} T' \xrightarrow{g} Y \to 0.
\]
Since $X \not\in {\rm add}(T)$, this sequence does not split.
Thus $X \not\cong Y$, because $X$ is rigid.
By Lemma \ref{prop1} we know that
$T \oplus Y$ is rigid.

Using Lemma \ref{approx2} and Lemma \ref{approx4}
yields that
$g$ is a minimal right ${\rm add}(T)$-approximation.
Thus, if $Y \in {\rm add}(T)$, then $T' \cong Y$ and $g$ would be an 
isomorphism.
But this would imply $X = 0$,
a contradiction since $X$ is indecomposable.
Thus $Y \not\in {\rm add}(T)$.

Next, we prove that
$Y$ is indecomposable.
Assume $Y = Y_1 \oplus Y_2$ with $Y_1$ and $Y_2$ non-zero.
For $i=1,2$ let $f_i\colon T_i' \to Y_i$ be a minimal right
${\rm add}(T)$-approximation, and let $X_i$ be the kernel.
Thus we have short exact sequences
\[
0 \to X_i \to T_i' \to Y_i \to 0.
\]
The direct sum $T_1' \oplus T_2' \to Y_1 \oplus Y_2$
is then a minimal right ${\rm add}(T)$-approximation.
Thus by the dual of Lemma \ref{approx3} the kernel $X_1 \oplus X_2$ 
is isomorphic to $X$.
Thus $X_1 = 0$ or $X_2 = 0$.
If $X_1 = 0$, then $0 \to T_1'$ is a direct summand
of $f\colon X \to T'$.
This is a contradiction to $f$ being minimal.
Similarly $X_2 = 0$ also leads to a contradiction.
Thus $Y$ must be indecomposable.
\end{proof}

Note that the assumption $\la \in {\rm add}(T)$ in
Proposition \ref{compl1} can be replaced by the weaker assumption
that there exists a monomorphism from $X$ to some object in 
${\rm add}(T)$.

The proof of
Proposition \ref{compl1} is similar to the proofs of
Lemma 6.3 - Lemma 6.6 in \cite{BMRRT}.
For convenience we gave a complete proof.
Note that we work with modules over preprojective
algebras, whereas \cite{BMRRT} deals with cluster categories.
However both have the crucial property that for all objects $M$ and $N$ the
extension groups $\Ext^1(M,N)$ and $\Ext^1(N,M)$ have the same dimension.

In the situation of the above proposition, we call $\{ X,Y \}$ an
{\it exchange pair associated to} $T$, and we write
\[
\mu_X(T \oplus X) = T \oplus Y.
\]
We say that $T \oplus Y$ is the mutation of $T \oplus X$ in direction $X$.
The short exact sequence
\[
0 \to X \xrightarrow{f} T' \xrightarrow{g} Y \to 0
\]
is the {\it exchange sequence} starting in $X$ and ending in $Y$.

For example, if $T = \la$ and $\{X,Y\}$ is an exchange pair associated 
to $T$, then 
\[
\mu_X(T \oplus X) = T \oplus \Omega^{-1}(X)
\]
where $\Omega$ is the syzygy functor.

Exchange sequences appear in tilting theory,
compare for example 
\cite[Theorem 1.1]{HU}, \cite[Section 3]{RS}, \cite[Proposition 1.3]{RS2} 
and \cite[Theorem 2.1]{U}.
A special case of an exchange sequence in the context of
tilting theory can be found already in \cite[Lemma 1.6]{APR}.

\begin{Prop}\label{compl2}
Let $X$ and $Y$ be indecomposable rigid $\la$-modules with 
$\dm \Ext_\la^1(Y,X) = 1$,
and let
\[
0 \to X \xrightarrow{f} M \xrightarrow{g} Y \to 0
\]
be a non-split short exact sequence.
Then $M \oplus X$ and $M \oplus Y$ are rigid
and $X,Y \notin {\rm add}(M)$.
If we assume additionally that $T \oplus X$ and
$T \oplus Y$ are basic maximal rigid $\la$-modules for some $T$,
then 
$f$ is a minimal left ${\rm add}(T)$-approximation
and
$g$ is a minimal right ${\rm add}(T)$-approximation.
\end{Prop}

Before we prove Proposition \ref{compl2} let us state a corollary:

\begin{Cor}\label{corcompl2}
Let $\{X,Y\}$ be an exchange pair associated to some basic rigid module
$T$ such that $T \oplus X$ and $T \oplus Y$ are maximal rigid, 
and assume $\dm \Ext_\la^1(Y,X) = 1$.
Then
\[
\mu_Y(\mu_X(T \oplus X)) = T \oplus X.
\]
\end{Cor}

\begin{proof}
Let 
\[
0 \to X \to T' \to Y \to 0
\]
be the short exact sequence from Proposition \ref{compl1}.
Thus $\mu_X(T \oplus X) = T \oplus Y$.
Since $\dm \Ext_\la^1(Y,X) = \dm \Ext_\la^1(X,Y) = 1$ 
and since $T \oplus X$ and $T \oplus Y$ are 
maximal rigid,
Proposition \ref{compl2} yields a non-split short exact sequence
\[
0 \to Y \xrightarrow{h} M \to X \to 0
\]
with $h$ a minimal left ${\rm add}(T)$-approximation.
Thus $\mu_Y(T \oplus Y) = T \oplus X$.
\end{proof}

\subsection{Proof of Proposition \ref{compl2}}
Let
$X$ and $Y$ be indecomposable rigid $\la$-modules with 
$\dm \Ext_\la^1(Y,X) = 1$,
and let
\begin{equation}\label{sequence}
0 \to X \xrightarrow{f} M \xrightarrow{g} Y \to 0
\end{equation}
be a non-split short exact sequence.

\begin{Lem}
$\Ext_\la^1(M,M) = 0$.
\end{Lem}

\begin{proof}
By Lemma \ref{extandhom} we have
\begin{align*}
\dm \Ext_\la^1(X \oplus Y,X \oplus Y) &= 
2\, \dm \Hom_\la(X \oplus Y,X \oplus Y) - 
(\dimv(X \oplus Y),\dimv(X \oplus Y)),\\
\dm \Ext_\la^1(M,M) &= 2\, \dm \Hom_\la(M,M) - (\dimv(M),\dimv(M)).
\end{align*}
Then our assumptions on $X$ and $Y$ yield
\[
2 = \dm \Ext_\la^1(X \oplus Y,X \oplus Y) =  2\, 
\dm \Hom_\la(X \oplus Y,X \oplus Y) - 
(\dimv(M),\dimv(M)).
\] 
Since Sequence~(\ref{sequence})
does not split, we get 
$M <_{\rm deg} X \oplus Y$, where $\leq_{\rm deg}$ is the usual
degeneration order, see for example \cite{Rie}.
Thus $\dm \Hom_\la(M,M) < \dm \Hom_\la(X \oplus Y,X \oplus Y)$, which implies
$\Ext_\la^1(M,M) = 0$.
\end{proof}

\begin{Lem}
$X,Y \notin {\rm add}(M)$. 
\end{Lem}

\begin{proof}
Assume $X \in {\rm add}(M)$.
Since $X$ is indecomposable, 
$M \cong X \oplus M'$ for some $M'$, and we get a
short exact sequence
\[
0 \to X \to X \oplus M' \to Y \to 0.
\]
By \cite[Proposition 3.4]{Rie} we get $M' \leq_{\rm deg} Y$.
Since $\Ext_\la^1(Y,Y) = 0$ this implies $M' = Y$.
Thus the above sequence splits, a contradiction. 
Dually, one shows that $Y \notin {\rm add}(M)$.
\end{proof}

\begin{Lem}
$\Ext_\la^1(M,X \oplus Y) = 0$.
\end{Lem}

\begin{proof}
Apply $\Hom_\la(-,X)$ to Sequence~(\ref{sequence}).
This yields an exact sequence
\begin{multline*}
0 \to \Hom_\la(Y,X) \to \Hom_\la(M,X) \xrightarrow{\Hom_\la(f,X)}
\Hom_\la(X,X) \xrightarrow{\delta} \Ext_\la^1(Y,X)\\ 
\to \Ext_\la^1(M,X) 
\to \Ext_\la^1(X,X) = 0.
\end{multline*}
Suppose that $\Hom_\la(f,X)$ is surjective. 
Then the identity morphism $X \to X$ factors through $f\colon X \to M$.
Thus $X \in {\rm add}(M)$, a contradiction to the previous lemma.
So the morphism $\delta$ has to be non-zero.
Since $\dm \Ext_\la^1(Y,X) = 1$ this implies that $\delta$ is surjective,
thus $\Ext_\la^1(M,X) = 0$.
Dually, one proves that $\Ext_\la^1(M,Y) = 0$.
\end{proof}

Now assume additionally that $T \oplus X$ and $T \oplus Y$ are 
basic maximal rigid for some $T$.

\begin{Lem}
$\Ext_\la^1(M,T) = 0$.
\end{Lem}

\begin{proof}
Applying $\Hom_\la(-,T)$ to Sequence~(\ref{sequence})
yields an exact sequence
\[
0 = \Ext_\la^1(Y,T) \to \Ext_\la^1(M,T) \to \Ext_\la^1(X,T) = 0.
\]
Thus $\Ext_\la^1(M,T) = 0$.
\end{proof}

\begin{Lem}
$M \in {\rm add}(T)$.
\end{Lem}

\begin{proof}
We know already that $X$ and $Y$ cannot be direct summands of $M$,
and what we proved so far yields that
$T \oplus X \oplus M$ is rigid.
Since $T \oplus X$ is a maximal rigid module, we get  
$M \in {\rm add}(T)$.
\end{proof}

Lemma \ref{approx2} and Lemma \ref{approx4} imply that
$f$ is a minimal left ${\rm add}(T)$-approximation.
Dually, $g$ is a minimal right ${\rm add}(T)$-approximation.
This finishes the proof of Proposition \ref{compl2}.


\section{Endomorphism algebras of maximal
rigid modules}\label{endosection}


In this section we work only with 
basic rigid $\la$-modules.
However, all our results on their endomorphism algebras are Morita
invariant, thus they 
hold for endomorphism algebras of arbitrary rigid $\la$-modules.

\subsection{Global dimension and quiver shapes}

\begin{Lem}\label{noloopcrit2}
Let $\{ X,Y \}$ be an exchange pair associated to a
basic rigid $\la$-module $T$.
Then the following are equivalent:
\begin{itemize}

\item
The quiver of $\End_\la(T \oplus X)$ has no loop at $X$;

\item
Every non-isomorphism $X \to X$ factors through ${\rm add}(T)$;

\item
$\dm \Ext_\la^1(Y,X) = 1$.

\end{itemize}
\end{Lem}

\begin{proof}
The equivalence of the first two statements is easy to show, we 
leave this to the reader.
Let 
\[
0 \to X \xrightarrow{f} T' \to Y \to 0
\]
be the exchange sequence starting in $X$.
Applying $\Hom_\la(-,X)$ yields an exact sequence
\[
0 \to \Hom_\la(Y,X) \to \Hom_\la(T',X) \xrightarrow{\Hom_\la(f,X)} 
\Hom_\la(X,X) \to \Ext_\la^1(Y,X) \to 0.
\]
Since $f$ is an ${\rm add}(T)$-approximation, every
non-isomorphism $X \to X$ factors through ${\rm add}(T)$ if and only
if it factors through $f$.
Clearly, this is equivalent to 
the cokernel $\Ext_\la^1(Y,X)$ of $\Hom_\la(f,X)$ being 1-dimensional.
Here we use that $\field$ is algebraically closed, which implies
$\Hom_\la(X,X)/{\rm rad}_\la(X,X) \cong \field$.
\end{proof}

\begin{Prop}\label{gldim}
Let $T$ be a basic maximal rigid $\la$-module.
If the quiver of $\End_\la(T)$ has no loops, then
\[
\gldim(\End_\la(T)) = 3.
\]
\end{Prop}

\begin{proof}
Set $E = \End_\la(T)$.
By assumption, the quiver of $E$ has no loops.
Thus $\Ext_E^1(S,S) = 0$ for all simple
$E$-modules $S$.
Let
\[
T = T_1 \oplus \cdots \oplus T_r
\]
with $T_i$ indecomposable for all $i$.
As before,
denote the simple $E$-module corresponding to $T_i$ by $S_{T_i}$.

Assume that $X = T_i$ is non-projective.
We claim that $\pdim_E(S_X) = 3$.
Let $\{X,Y\}$ be the exchange pair associated to $T/X$.
Note that $\la \in {\rm add}(T/X)$.
By Lemma \ref{noloopcrit2} we have $\dm \Ext_\la^1(Y,X) = 1$.
Let
\[
0 \to X \xrightarrow{f} T' \to Y \to 0
\]
and
\[
0 \to Y \to T'' \to X \to 0
\]
be the corresponding non-split short exact sequences.
Applying $\Hom_\la(-,T)$ to both sequences yields  
exact sequences
\begin{multline*}
0 \to \Hom_\la(Y,T) \to \Hom_\la(T',T) \xrightarrow{\Hom_\la(f,T)} 
\Hom_\la(X,T)\\ 
\to \Ext_\la^1(Y,T) \cong \Ext_\la^1(Y,X) \to \Ext_\la^1(T',T) = 0
\end{multline*}
and
\[
0 \to \Hom_\la(X,T) \to \Hom_\la(T'',T) \to \Hom_\la(Y,T) 
\to \Ext_\la^1(X,T) = 0.
\]
Since $\dm \Ext_\la^1(Y,X) = 1$ we know that the cokernel of 
$\Hom_\la(f,T)$ is
1-dimensional.
Thus, since $\Hom_\la(X,T)$ is the indecomposable projective $E$-module with
top $S_X$, the cokernel must be isomorphic to $S_X$.
Combining the two sequences above yields an exact sequence
\[
0 \to \Hom_\la(X,T) \to \Hom_\la(T'',T) \to \Hom_\la(T',T) \to \Hom_\la(X,T) 
\to S_X \to 0. 
\]
This is a projective resolution of $S_X$.
Thus $\pdim_E(S_X) \le 3$.

By Proposition \ref{compl1} we have $X \notin {\rm add}(T'')$.
Thus $\Hom_E(\Hom_\la(T'',T),S_X) = 0$ and
$\Ext_E^3(S_X,S_X) \cong \Hom_E(\Hom_\la(X,T),S_X)$
is one-dimensional, in particular it is non-zero.
Thus $\pdim_E(S_X) = 3$.

Next, assume that $P = T_i$ is projective.
We claim that $\pdim_E(S_P) \le 2$.
Set 
$X = P/S$ where $S$ is the (simple) socle of $P$.
First, we prove that $X$ is rigid:
Applying $\Hom_\la(-,X)$ to the short exact sequence
\[
0 \to S \to P \xrightarrow{\pi} X \to 0
\]
yields an exact sequence
\[
0 \to \Hom_\la(X,X) \to \Hom_\la(P,X) \to \Hom_\la(S,X) 
\to \Ext_\la^1(X,X) \to \Ext_\la^1(P,X) = 0.
\]
The quiver of the preprojective algebra $\la$ does not contain
any loops.
Thus the socle of $X$ does not contain $S$ as a composition factor.
This implies $\Hom_\la(S,X) = 0$, and thus $\Ext_\la^1(X,X) = 0$.
Let $f\colon X \to T'$ be a minimal left ${\rm add}(T/P)$-approximation.
Clearly, $f$ is injective. 
We get a short exact sequence
\[
0 \to X \xrightarrow{f} T' \to Y \to 0.
\]
Now Lemma \ref{prop1} yields that
$Y \oplus T/P$ is rigid.
Since $P$ is projective, there is only one indecomposable
module $C$ such that $C \oplus T/P$ is maximal rigid, namely $C = P$.
Here we use the assumption that $T$ is maximal rigid.
This implies
$Y \in {\rm add}(T)$.
The projection $\pi\colon P \to X$ yields an exact sequence
\[
P \xrightarrow{h} T' \to Y \to 0
\]
where $h = f\pi$.
Applying $\Hom_\la(-,T)$ to this sequence gives an exact sequence
\[
0 \to \Hom_\la(Y,T) \to \Hom_\la(T',T) \xrightarrow{\Hom_\la(h,T)} 
\Hom_\la(P,T) \to Z \to 0.
\]
We have $\Hom_\la(P,T) = \Hom_\la(P,T/P) \oplus \Hom_\la(P,P)$.
For each morphism $g\colon P \to T/P$ there exists some
morphism $g'\colon X \to T/P$ such that $g = g'\pi$.
Since $f$ is a left ${\rm add}(T/P)$-approximation of $X$, there
is a morphism $g''\colon T' \to T/P$ such that $g' = g''f$.
Thus we get a commutative diagram
\[
\xymatrix{
P \ar[r]^{g} \ar[d]_{\pi}    & T/P \\
X \ar[ur]^{g'} \ar[r]_{f}    & T' \ar[u]_{g''}.
}
\]
This implies $g = g'\pi = g''f\pi = g''h$.
Thus $g$ factors through $h$.

Since there is no loop at $S_P$, all non-isomorphisms 
$P \to P$ factor through $T/P$.
Thus by the argument above they factor through $h$.
Thus the cokernel $Z$ of $\Hom_\la(h,T)$ must be 1-dimensional, which
implies $Z \cong S_P$.
Since $Y \in {\rm add}(T)$, we know that $\Hom_\la(Y,T)$ is projective.
Thus the above is a projective resolution of $S_P$, 
so $\pdim_E(S_P) \le 2$.

This finishes the proof of Proposition \ref{gldim}.
\end{proof}

\begin{Cor}
$\repdim(\la) \le 3$.
\end{Cor}

\begin{proof}
Clearly, the module $T_Q$ from Section \ref{deftdelta} 
is a generator-cogenerator of $\md(\la)$.
The quiver of $T_Q$ has no loops, thus
by Proposition \ref{gldim} we know that 
$\gldim(\End_\la(T_Q)) = 3$.
This implies $\repdim(\la) \le 3$.
\end{proof}

The statements in the following theorem are presented in the order
in which we prove them.

\begin{Thm}\label{quivershape}
Let $T$ be a basic maximal rigid $\la$-module, and set $E = \End_\la(T)$.
Then the following hold:
\begin{itemize}

\item[(1)]
The quiver of $E$ has no loops;

\item[(2)]
$\gldim(E) = 3$;

\item[(3)]
$T$ is maximal 1-orthogonal;

\item[(4)]
$\domdim(E) = 3$;

\item[(5)]
The quiver of $E$ has no sinks and no sources;

\item[(6)]
For all simple $E$-modules $S$ we have
$\Ext_E^1(S,S) = 0$ and
$\Ext_E^2(S,S) = 0$;

\item[(7)]
The quiver of $E$ has no 2-cycles.

\end{itemize}
\end{Thm}

\begin{proof}
By Theorem \ref{tilting} we know that $\End_\la(T_Q)$ and
$\End_\la(T)$ are derived equivalent, where $T_Q$ is the
complete rigid module mentioned in Section \ref{deftdelta}.
Since the quiver of $\End_\la(T_Q)$ has no loops, 
Proposition \ref{gldim} implies that $\gldim(\End_\la(T_Q)) = 3 < \infty$.
This implies
$\gldim(\End_\la(T)) < \infty$.
Thus by Theorem \ref{Igusa} the quiver of $\End_\la(T)$ has no loops.
Then again Proposition \ref{gldim} yields $\gldim(\End_\la(T)) = 3$.
Thus $T$ is maximal 1-orthogonal by Theorem \ref{Iyama1}, and
by Theorem \ref{Iyama2} we get $\domdim(\End_\la(T)) = 3$.
(It follows from the definitions that for an algebra $A$ and some $n \ge 1$,
$\gldim(A) = n$ implies 
$\domdim(A) \le n$.)
This proves parts (1)-(4) of the theorem.

For any indecomposable direct summand $M$ of $T$ there
are non-zero homomorphisms $M \to P_i$ and
$P_j \to M$ for some indecomposable projective $\la$-modules
$P_i$ and $P_j$ which are not isomorphic to $M$.
Thus the vertex $S_M$ in the quiver of $E$  
is neither a sink nor a source.
So (5) is proved.

Since the quiver of $E$ has no loops, we have $\Ext_E^1(S,S) = 0$
for all simple $E$-modules $S$.
Let $X$ be a non-projective direct summand of $T$.
In the proof of Proposition \ref{gldim}, we constructed a projective resolution
\[
0 \to \Hom_\la(X,T) \to \Hom_\la(T'',T) \to \Hom_\la(T',T) 
\to \Hom_\la(X,T) \to S_X \to 0, 
\]
and we also know that $X \notin {\rm add}(T'')$.
Thus applying $\Hom_E(-,S_X)$ to this resolution yields
$\Ext_E^2(S_X,S_X) = 0$.
Next, assume
$P$ is an indecomposable projective direct summand of $T$. 
As in the proof of Proposition \ref{gldim} we have a projective resolution
\[
0 \to \Hom_\la(Y,T) \to \Hom_\la(T',T) \xrightarrow{\Hom_\la(h,T)} 
\Hom_\la(P,T) \to S_P \to 0
\]
where $P \notin {\rm add}(T')$.
Since the module $T'$ projects onto $Y$, we conclude that
$P \notin {\rm add}(Y)$.
Applying $\Hom_E(-,S_P)$ to the above resolution of $S_P$ yields
$\Ext_E^2(S_P,S_P) = 0$.
This finishes the proof of (6).

We proved that
$\Ext_E^2(S,S) = 0$ for all simple $E$-modules $S$.
We also know that $\gldim(E) = 3 < \infty$.
Then it follows from Proposition \ref{2cycleresult} that 
the quiver of $E$ cannot have 2-cycles.
Thus (7) holds. 
This finishes the proof.
\end{proof}

\begin{Cor}\label{corquivershape}
Let $T = T_1 \oplus \cdots \oplus T_r$ be a basic maximal rigid
$\la$-module with $T_i$ indecomposable
for all $i$.
For a non-projective $X = T_i$ let 
\[
0 \to X \to T' \to Y \to 0
\]
be the corresponding exchange sequence starting in $X$.
Then the following hold:
\begin{itemize}

\item
We have $\dm \Ext_\la^1(Y,X) = \dm \Ext_la^1(X,Y) = 1$, and 
the exchange sequence ending in $X$ is of the form
\[
0 \to Y \to T'' \to X \to 0
\]
for some $T'' \in {\rm add}(T/X)$;

\item
The simple $\End_\la(T)$-module $S_X$ has a minimal projective resolution
of the form
\[
0 \to \Hom_\la(X,T) \to \Hom_\la(T'',T) \to \Hom_\la(T',T) \to \Hom_\la(X,T) 
\to S_X \to 0;
\]

\item
We have ${\rm add}(T') \cap {\rm add}(T'') = 0$.

\end{itemize}
\end{Cor}

\begin{proof}
By Theorem \ref{quivershape}
the quiver of $\End_\la(T)$
has no loops.
Now Lemma \ref{noloopcrit2} yields that $\dm \Ext_\la^1(Y,X) = 1$, and by
Lemma \ref{extandhom} we get $\dm \Ext_\la^1(X,Y) = 1$.
Corollary \ref{corcompl2} implies that the 
exchange sequence
ending in $X$ starts at $Y$.
The minimal projective resolution of $S_X$
is obtained from the proof of Proposition \ref{gldim}.
By Theorem \ref{quivershape} there are no 2-cycles in the
quiver of $\End_\la(T)$.
This implies ${\rm add}(T') \cap {\rm add}(T'') = 0$.
\end{proof}

\subsection{Ext-group symmetries}

\begin{Prop}\label{CalabiYau}
Let $T$ be a basic maximal rigid $\la$-module, and let $X$ be
a non-projective indecomposable direct summand of $T$.
Set $E = \End_\la(T)$.
Then for any simple $E$-module $S$ we have
\[
\dm \Ext_E^{3-i}(S_X,S) = \dm \Ext_E^i(S,S_X)
\]
where $0 \le i \le 3$.
\end{Prop}

\begin{proof}
We have $S = S_Z$ for some indecomposable direct summand $Z$ of $T$.
Since $S_X$ and $S_Z$ are simple, we get 
\[
\dm \Hom_E(S_X,S_Z) = \dm \Hom_E(S_Z,S_X) =
\begin{cases}
1 & \text{if $S_Z \cong S_X$},\\
0 & \text{otherwise}.
\end{cases}
\]
Let
\begin{equation}\label{projres}
0 \to \Hom_\la(X,T) \to \Hom_\la(T'',T) \to \Hom_\la(T',T) 
\to \Hom_\la(X,T) \to S_X \to 0
\end{equation}
be the minimal projective resolution of $S_X$ as constructed in the proof
of Proposition \ref{gldim}.
Dually we get a minimal injective resolution
\begin{equation}\label{injres}
0 \to S_X \to {\rm D}\Hom_\la(T,X) \to {\rm D}\Hom_\la(T,T'') \to 
{\rm D}\Hom_\la(T,T') \to {\rm D}\Hom_\la(T,X) \to 0
\end{equation}
of $S_X$.
We apply $\Hom_E(-,S_Z)$ to (\ref{projres}) and 
$\Hom_E(S_Z,-)$ to (\ref{injres})
and get
\[
\dm \Ext_E^3(S_X,S_Z) = \dm \Ext_E^3(S_Z,S_X) =
\begin{cases}
1 & \text{if $S_Z \cong S_X$},\\
0 & \text{otherwise}.
\end{cases}
\]
Here we use that $X \notin {\rm add}(T' \oplus T'')$.
Similarly, if $S_Z \cong S_X$, then 
$\Ext_E^1(S_X,S_Z) = 0 = \Ext_E^2(S_X,S_Z)$ and 
$\Ext_E^1(S_Z,S_X) = 0 = \Ext_E^2(S_Z,S_X)$.
Thus the proposition is true for $S_Z \cong S_X$.

From now on assume $S_X \not\cong S_Z$.
The projective resolution (\ref{projres}) yields a complex
\[
0 \to \Hom_E(\Hom_\la(T',T),S_Z) \xrightarrow{\theta} 
\Hom_E(\Hom_\la(T'',T),S_Z) \to 0
\]
whose homology groups are the extension groups $\Ext_E^i(S_X,S_Z)$.
By Corollary \ref{corquivershape} we have 
${\rm add}(T') \cap {\rm add}(T'') = 0$, which implies that
either $\Hom_E(\Hom_\la(T',T),S_Z) = 0$ or 
$\Hom_E(\Hom_\la(T'',T),S_Z) = 0$.
Hence
$\theta = 0$.
Therefore we have $\dm \Ext_E^2(S_X,S_Z) = [T'':Z]$.
The discussion in Section \ref{endoquiver} yields
\begin{align*}
\dm \Ext_E^1(S_X,S_Z) &= [T':Z] = \text{number of arrows $S_X \to S_Z$ in 
$\Gamma_T$},\\
\dm \Ext_E^1(S_Z,S_X) &= [T'':Z] = \text{number of arrows $S_Z \to S_X$ 
in $\Gamma_T$}.
\end{align*}
(As before $\Gamma_T$ denotes the quiver of $E = \End_\la(T)$.)

This implies
$\dm \Ext_E^2(S_X,S_Z) = \dm \Ext_E^1(S_Z,S_X)$.
Using the injective resolution (\ref{injres}) we get
$\dm \Ext_E^2(S_Z,S_X) = \dm \Ext_E^1(S_X,S_Z)$.
This finishes the proof.
\end{proof}

\subsection{The graph of basic maximal rigid modules}

\begin{Prop}\label{numberofcompl}
Let $T$ be a basic maximal rigid $\la$-module, and let
$X$ be an indecomposable direct summand of $T$.
If $X$ is non-projective, then up to isomorphism there
exists exactly one
indecomposable $\la$-module 
$Y$ such that $X \not\cong Y$ and $Y \oplus T/X$ is maximal rigid.
\end{Prop}

\begin{proof}
Recall that every maximal rigid $\la$-module
is maximal 1-orthogonal.
Proposition \ref{compl1} yields an exchange sequence
\[
0 \to X \to T' \to Y \to 0
\]
such that $Y$ satisfies the required properties.
Proposition \ref{numberofcomplements1} implies that $Y$ is uniquely determined
up to isomorphism.
\end{proof}

Let $T$ be a maximal rigid $\la$-module.
Let ${\mathcal T}_{\End_\la(T)}$ be the graph of basic 
tilting modules over $\End_\la(T)$.
By definition the vertices of this graph are the isomorphism classes
of basic tilting modules over $\End_\la(T)$.
Two such basic tilting modules $M_1$ and $M_2$ are connected with an
edge if and only if $M_1 = M \oplus M_1'$ and $M_2 = M \oplus M_2'$
for some $M$ and some indecomposable
modules $M_1'$ and $M_2'$ with $M_1' \not\cong M_2'$.

Similarly, let ${\mathcal T}_\la$ be the graph with set of vertices
the isomorphism classes of basic maximal rigid $\la$-modules,
and an edge between vertices $T_1$ and $T_2$ if and only if
$T_1 = T \oplus T_1'$ and $T_2 = T \oplus T_2'$ for some $T$ and
some indecomposable modules $T_1'$ and $T_2'$ with $T_1' \not\cong T_2'$.

\begin{Lem}\label{tiltgraph1}
Let $T =  T_1 \oplus \cdots \oplus T_r$ be a basic maximal rigid
$\la$-module with $T_i$ indecomposable for all $i$, and 
assume that $T_{r-n+1},\ldots,T_r$ are projective.
There are exactly $n$ indecomposable projective-injective
$\End_\la(T)$-modules up to isomorphism, namely 
$
\Hom_\la(T_{r-n+1},T),\ldots,\Hom_\la(T_r,T).
$
These are direct summands of any tilting module over $\End_\la(T)$.
\end{Lem}

\begin{proof}
For $1 \le i \le r-n$ 
we know that
$\Hom_\la(\mu_{T_i}(T),T)$ is a tilting module over $\End_\la(T)$, which 
does not contain $\Hom_\la(T_i,T)$ as a direct summand.
Thus $\Hom_\la(T_i,T)$ cannot be projective-injective,
otherwise it would occur as a direct summand of
any tilting module over $\End_\la(T)$.

Let $\nu_\la = {\rm D}\Hom_\la(-,\la)$ be the Nakayama automorphism of
$\md(\la)$.
It is well known that $\nu_\la$ maps projective modules to injective ones.
But $\la$ is a selfinjective algebra.
This implies that for $r-n+1 \le i \le r$ we get
\[
\Hom_\la(T_i,T) \cong {\rm D}\Hom_\la(T,\nu_\la(T_i)) = 
{\rm D}\Hom_\la(T,T_{j})
\]
for some $j$.
In particular, $\Hom_\la(T_i,T)$ is a projective-injective 
$\End_\la(T)$-module.
\end{proof}

\begin{Prop}\label{tiltgraph2}
Let $T$ be a basic maximal rigid $\la$-module.
The functor 
\[
F_T\colon \md(\la) \to \md(\End_\la(T))
\]
induces an embedding of graphs 
\[
\psi_T\colon {\mathcal T}_\la \to {\mathcal T}_{\End_\la(T)}
\]
whose image is a union of connected components
of ${\mathcal T}_{\End_\la(T)}$.
Each vertex of ${\mathcal T}_\la$ (and therefore each vertex of
the image of $\psi_T$) has exactly $r-n$ neighbours.
\end{Prop}

\begin{proof}
By Theorem \ref{Iyama3} we know that every vertex of ${\mathcal T}_\la$
gets mapped to a vertex of ${\mathcal T}_{\End_\la(T)}$.
Proposition \ref{homfunctor} implies that $F_T$ 
reflects isomorphism classes, therefore $\psi_T$
is injective on vertices.
Proposition \ref{compl1} and
Proposition \ref{numberofcompl} and its proof yield that $\psi_T$ is
injective on edges as well.
It also follows that every vertex of ${\mathcal T}_\la$
has exactly $r-n$ neighbours.
Thus every vertex in the image of $\psi_T$ has at least $r-n$ neighbours.
But by \cite[Proposition 1.3]{RS2} there are at most two complements
for an almost complete tilting module.
Combining this with Lemma \ref{tiltgraph1} implies that
every vertex in the image of $\psi_T$ has exactly 
$r-n$ neighbours.
Thus the image of $\psi_T$ is a union of connected components.
\end{proof}

\begin{Conj}\label{connected}
The graphs ${\mathcal T}_\la$ and ${\mathcal T}_{\End_\la(T)}$ are connected.
\end{Conj}

To prove Conjecture \ref{connected} it is enough to show that 
${\mathcal T}_{\End_\la(T)}$ is connected.


\section{From mutation of modules to 
mutation of matrices}\label{graphmutation}


In this section we prove Theorem \ref{Mainresult1}.

\subsection{}
Let $B = (b_{ij})$ be
an $l \times m$-matrix with real entries such that $l \ge m$,
and let $k \in [1,m]$.
Following Fomin and Zelevinsky,
the {\it mutation} of $B$ in direction $k$ is an $l \times m$-matrix 
\[
\mu_k(B) = (b_{ij}')
\]
defined by
\[
b_{ij}' =
\begin{cases}
-b_{ij} & \text{if $i=k$ or $j=k$},\\
b_{ij} + \dfrac{|b_{ik}|b_{kj} + b_{ik}|b_{kj}|}{2} & \text{otherwise},
\end{cases}
\]
where $i \in [1,l]$ and $j \in [1,m]$.
If $l = m$ and $B$ is skew-symmetric, then 
it is easy to check that $\mu_k(B)$ is also a skew-symmetric 
matrix.

For an $m \times m$-matrix $B$ and some
$k \in [1,m]$ we define an $m \times m$-matrix $S = S(B,k) = (s_{ij})$
by
\[
s_{ij} =
\begin{cases}
-\delta_{ij} + \dfrac{|b_{ij}|-b_{ij}}{2} & \text{if $i = k$},\\
\delta_{ij} & \text{otherwise}.
\end{cases}
\]

By $S^t$ we denote the transpose of the matrix $S = S(B,k)$.
If $B$ is skew-symmetric, then one easily checks that 
$$S^2 = 1$$ 
is the 
identity matrix.
The following lemma is a special case of \cite[(3.2)]{BFZ}.
For convenience we include a proof.

\begin{Lem}\label{mutation1}
Let $B$ be a skew-symmetric $m \times m$-matrix, and let
$S = S(B,k)$ for some $k \in [1,m]$.
Then we have
\[
\mu_k(B) = S^tBS.
\]
\end{Lem}

\begin{proof}
Let $b_{ij}'$ be the $ij$th entry of $S^tBS$.
Thus
\[
b_{ij}' = \sum_{p=1}^m \sum_{q=1}^m s_{pi} b_{pq} s_{qj}.
\]
Since $B$ is skew-symmetric we have $b_{ll} = 0$ for all $l$.
From the definition of $S$ we get 
\[
b_{ij}' = 
\begin{cases}
s_{kk} b_{kj} s_{jj} & \text{if $i = k$},\\
s_{ii} b_{ik} s_{kk} & \text{if $j = k$},\\
s_{ii} b_{ij} s_{jj} + s_{ii} b_{ik} s_{kj} + s_{ki} b_{kj} s_{jj} 
& \text{otherwise}.
\end{cases}
\] 
Now an easy calculation yields $\mu_k(B) = S^tBS$.
\end{proof}

\subsection{}
Let $A$ be a finite-dimensional algebra, and let
$P_1, \ldots, P_m$ be a complete set of representatives of
isomorphism classes of indecomposable projective $A$-modules.
By $S_i$ we denote the top of $P_i$.
Thus $S_1, \ldots, S_m$ is a complete set of representatives of
isomorphism classes of simple $A$-modules.
The {\it Cartan matrix} $C$ of $A$ is by definition
$C = (c_{ij})_{1 \le i,j \le m}$ where
\[
c_{ij} = \dm \Hom_A(P_i,P_j).
\]
Now assume that the global dimension of $A$ is finite.
Then the {\it Ringel form}
of $A$ is a bilinear form $\bil{-,-}\colon \Z^m \times \Z^m \to \Z$
defined by
\[
\bil{\dimv(M),\dimv(N)} = \sum_{i \ge 0} (-1)^i\, \dm \Ext_A^i(M,N)
\]
where $M, N \in \md(A)$.
Here we use the convention $\Ext_A^0(M,N) = \Hom_A(M,N)$.
We often just write $\bil{M,N}$ instead of $\bil{\dimv(M),\dimv(N)}$.
The matrix of the Ringel form of $A$ is then
\[
R = (\bil{S_i,S_j})_{1 \le i,j \le m}.
\]
The following lemma can be found in \cite[Section 2.4]{Ri1}:

\begin{Lem}\label{ringelform}
If $\gldim(A) < \infty$, then the Cartan matrix $C$ of $A$ is
invertible, and we have
\[
R = C^{-t}, 
\]
where $C^{-t}$ denotes the transpose of the inverse of $C$.
\end{Lem}

\subsection{}
Let $T = T_1 \oplus \cdots \oplus T_r$ be a basic complete rigid
$\la$-module with $T_i$ indecomposable for all $i$.
Without loss of generality we assume that 
$T_{r-n+1}, \ldots, T_r$ are projective.
For $1 \le i \le r$ let $S_i$ be the simple $\End_\la(T)$-module
corresponding to $T_i$.
The matrix
\[
C_T = (c_{ij})_{1 \le i,j \le r}
\]
where 
\[
c_{ij} = \dm \Hom_{\End_\la(T)}(\Hom_\la(T_i,T),\Hom_\la(T_j,T)) 
= \dm \Hom_\la(T_j,T_i)
\]
is the Cartan matrix of the algebra $\End_\la(T)$.

By Theorem \ref{quivershape} we know that $\gldim(\End_\la(T)) = 3$.
Thus by Lemma \ref{ringelform}
\[
R_T =(r_{ij})_{1 \le i,j \le r} = C_T^{-t}
\]
is the matrix of the Ringel form of $\End_\la(T)$, where 
\[
r_{ij} = \bil{S_i,S_j} = \sum_{i=0}^3 (-1)^i \, 
\dm \Ext_{\End_\la(T)}^i(S_i,S_j). 
\]

\begin{Lem}\label{rijequation}
Assume that $i \le r-n$ or $j \le r-n$.
Then the following hold:
\begin{itemize} 

\item
$r_{ij} = \dm \Ext_{\End_\la(T)}^1(S_j,S_i) - 
\dm \Ext_{\End_\la(T)}^1(S_i,S_j)$;

\item
$r_{ij} = - r_{ji}$;

\item
$
r_{ij} =
\begin{cases}
\text{number of arrows $j \to i$ in $\Gamma_T$}    & 
\text{if $r_{ij} > 0$},\\
-(\text{number of arrows $i \to j$ in $\Gamma_T$}) & 
\text{if $r_{ij} < 0$},\\
0                                   & \text{otherwise}.
\end{cases}
$

\end{itemize}
\end{Lem}

\begin{proof}
The first two statements follow from
Proposition \ref{CalabiYau}.
Since by Theorem \ref{quivershape} there are no 2-cycles in the quiver 
of $\End_\la(T)$, 
at least one of the two summands in the equation 
$$
r_{ij} = \dm \Ext_{\End_\la(T)}^1(S_j,S_i) - 
\dm \Ext_{\End_\la(T)}^1(S_i,S_j)
$$ 
has to be 0.
Now the third statement follows from the remarks in Section \ref{endoquiver}.
\end{proof}

Recall that $B(T) = (t_{ij})_{1 \le i,j \le r}$ is the $r \times r$-matrix 
defined by
\[
t_{ij} = 
(\text{number of arrows $j \to i$ in $\Gamma_T$}) 
-(\text{number of arrows $i \to j$ in $\Gamma_T$}).
\]
Let 
$B(T)^\circ = (t_{ij})$ and $R_T^\circ = (r_{ij})$
be the $r \times (r-n)$-matrices obtained from
$B(T)$ and $R_T$, respectively, by deleting the last $n$ columns.
As a consequence of Lemma \ref{rijequation} we get the following:

\begin{Cor}\label{corrijequation}
$R_T^\circ = B(T)^\circ$.
\end{Cor}

Note that 
the dimension vector of the indecomposable projective 
$\End_\la(T)$-module
$\Hom_\la(T_i,T)$ is the $i$th column of the matrix $C_T$.

For $1 \le k \le r-n$ let
\begin{equation}\label{seq1}
0 \to T_k \to T' \to T_k^* \to 0
\end{equation}
and
\begin{equation}\label{seq2}
0 \to T_k^* \to T'' \to T_k \to 0
\end{equation}
be exchange sequences associated to the direct summand $T_k$
of $T$.
Keeping in mind the remarks in Section \ref{endoquiver},
it follows from Lemma \ref{rijequation} that
\[
T'  = \bigoplus_{r_{ik} > 0} T_i^{r_{ik}}
\;\;\;\; \text{ and }  \;\;\;\;
T'' = \bigoplus_{r_{ik} < 0} T_i^{-r_{ik}}.
\]
Set 
\[
T^* = \mu_{T_k}(T) = T_k^* \oplus T/T_k
\]
and $S = S(R_T,k)$.

\begin{Prop}\label{mutation3}
With the above notation we have
\[
C_{T^*} = SC_TS^t.
\]
\end{Prop}

\begin{proof}
Set $C_T' = SC_TS^t = (c_{ij}')_{1 \le i,j \le r}$.
Thus we have
\[
c_{ij}' = \sum_{p=1}^r \sum_{q=1}^r s_{ip} c_{pq} s_{jq}
\]
where
\[
s_{ip} =
\begin{cases}
-\delta_{ip} + \dfrac{|r_{ip}|-r_{ip}}{2} & \text{if $i = k$},\\
\delta_{ip} & \text{otherwise}
\end{cases}
\]
and
\[
s_{jq} =
\begin{cases}
-\delta_{qj} + \dfrac{|r_{qj}|+r_{qj}}{2} & \text{if $j = k$},\\
\delta_{qj} & \text{otherwise}.
\end{cases}
\]
It follows that the transformation
$C_T \mapsto C_T'$ only changes the $k$th row and the $k$th column
of $C_T$.
We denote the $ij$th entry of $C_{T^*}$ by $c_{ij}^*$.
Clearly, $c_{ij}^* = c_{ij} = c_{ij}'$ 
provided $i \not= k$ and $j \not= k$.
Assume $i \not= k$ and $j = k$.
We get 
\[
c_{ik}' = \sum_{q=1}^r c_{iq}s_{kq} = 
- c_{ik} + \sum_{r_{qk}>0} r_{qk}c_{iq}.
\]
Applying $\Hom_\la(-,T_i)$ to Sequence~(\ref{seq1}) yields 
\begin{align*}
c_{ik}^* &= \dm \Hom_\la(T_k^*,T_i)
= - \dm \Hom_\la(T_k,T_i)
+ \sum_{r_{qk}>0} r_{qk}\, \dm \Hom_\la(T_q,T_i)\\
&= -c_{ik} + \sum_{r_{qk}>0} r_{qk}c_{iq}.
\end{align*}
Thus $c_{ik}' = c_{ik}^*$.
Similarly, if $i = k$ and $j \not= k$, then
\[
c_{kj}' = \sum_{p=1}^r s_{kp}c_{pj} = 
- c_{kj} + \sum_{r_{kp}<0} |r_{kp}|c_{pj}.
\]
Applying $\Hom_\la(T_j,-)$ to Sequence~(\ref{seq1}) yields
\begin{align*}
c_{kj}^* &= \dm \Hom_\la(T_j,T_k^*)
= - \dm \Hom_\la(T_j,T_k) + \sum_{r_{pk}>0} r_{pk}\, \dm \Hom_\la(T_j,T_p)\\
&= -c_{kj} + \sum_{r_{pk}>0} r_{pk}c_{pj}
= -c_{kj} + \sum_{r_{kp}<0} |r_{kp}|c_{pj}.
\end{align*}
Thus $c_{kj}' = c_{kj}^*$.
Finally, let $i = j = k$.
Thus
$$
c_{kk}' = \sum_{p=1}^r \sum_{q=1}^r s_{kp} c_{pq} s_{kq}
        = c_{kk} - \sum_{r_{kp}<0} |r_{kp}|c_{pk}    
- \sum_{r_{qk}>0} r_{qk}c_{kq} 
+ \sum_{r_{qk}>0} \sum_{r_{kp}<0} r_{qk}|r_{kp}|c_{pq}.
$$
We apply $\Hom_\la(-,T_k^*)$ to Sequence~(\ref{seq1}) and get
\[
c_{kk}^* = \dm \Hom_\la(T_k^*,T_k^*) = - \dm \Hom_\la(T_k,T_k^*) 
+ \sum_{r_{qk}>0} r_{qk}\, \dm \Hom_\la(T_q,T_k^*).
\]
By applying $\Hom_\la(T_k,-)$ and $\Hom_\la(T_q,-)$ to the same sequence, 
we can compute
the values $\dm \Hom_\la(T_k,T_k^*)$ and $\dm \Hom_\la(T_q,T_k^*)$ 
in the above equation.
We get
\begin{align*}
c_{kk}^* &= c_{kk} - \sum_{r_{pk}>0} r_{pk}c_{pk}  
+ \left( \sum_{r_{qk}>0} r_{qk} \left(-c_{kq} + \sum_{r_{pk}>0}r_{pk}c_{pq}  
\right) \right)\\
         &= c_{kk} - \sum_{r_{pk}>0} r_{pk}c_{pk}  
- \sum_{r_{qk}>0} r_{qk}c_{kq}  + 
\sum_{r_{qk}>0} \sum_{r_{pk}>0} r_{qk}r_{pk}c_{pq} \\ 
         &= c_{kk} - \sum_{r_{kp}<0} |r_{kp}|c_{pk}  
- \sum_{r_{qk}>0} r_{qk}c_{kq}  + 
\sum_{r_{qk}>0} \sum_{r_{kp}<0} r_{qk}|r_{kp}|c_{pq}.  
\end{align*}
This proves that $c_{kk}' = c_{kk}^*$.
\end{proof}

Note that in the proof of Proposition \ref{mutation3} we only made
use of Sequence~(\ref{seq1}).
There is an alternative proof using Sequence~(\ref{seq2}).

\begin{Cor}\label{mutation4}
$R_{T^*} = S^tR_TS$.
\end{Cor}

\begin{proof}
From Proposition \ref{mutation3} and 
and Lemma \ref{ringelform}
we get
\[
C_{T^*}^{-1} = (S^t)^{-1}C_T^{-1}S^{-1} = S^tC_T^{-1}S
\]
and therefore
\[
R_{T^*} = C_{T^*}^{-t} = S^tC_T^{-t}S = S^tR_TS.
\]
\end{proof}

\begin{Cor}\label{mutation5}
$R_{T^*}^\circ = \mu_k(R_T^\circ)$.
\end{Cor}

\begin{proof}
For an $r \times r$-matrix $B = (b_{ij})_{1 \le i,j \le r}$
define a matrix $B^\vee = (b_{ij}^\vee)_{1 \le i,j \le r}$ by
\[
b_{ij}^\vee = 
\begin{cases}
0      & \text{if $r-n+1 \le i,j \le r$},\\
b_{ij} & \text{otherwise}.
\end{cases}
\]
It follows from Lemma \ref{rijequation} that the matrix $R_T^\vee$
is skew-symmetric.
One easily checks that $S(R_T^\vee,k) = S(R_T,k) = S$.
By Corollary \ref{mutation4} we have $R_{T^*} = S^tR_TS$.

It follows from the definition that the matrix $S$ is of the
form
\[
S = \left( \begin{matrix}S_1 & S_2\\ {\bf 0} & {\bf I}  \end{matrix} \right)
\]
where ${\bf 0}$ is the 0-matrix of size $n \times (r-n)$,
and ${\bf I}$ is the identity matrix of size $n \times n$.

We partition $R_T$ and $R_T^\vee$ into blocks of the corresponding sizes
and obtain
\[
R_T = \left( \begin{matrix}R_1 & R_2\\R_3 & R_4  \end{matrix} \right)
\text{ and }
R_T^\vee = \left( \begin{matrix}R_1 & R_2\\R_3 & 0  \end{matrix} \right).
\]
This implies
\[
R_{T^*} = S^tR_TS = 
\left( \begin{matrix}
S_1^tR_1S_1           & S_1^tR_1S_2 + S_1^tR_2\\
S_2^tR_1S_1 + R_3S_1  & S_2^tR_1S_2 + S_2^tR_2 + R_3S_2 + R_4
\end{matrix} \right).
\]
Since $R_T^\vee$ is skew-symmetric and $S = S(R_T^\vee,k)$, 
Lemma \ref{mutation1} implies 
\[
\mu_k(R_T^\vee) = S^tR_T^\vee S =
\left( \begin{matrix}
S_1^tR_1S_1           & S_1^tR_1S_2 + S_1^tR_2\\
S_2^tR_1S_1 + R_3S_1  & S_2^tR_1S_2 + S_2^tR_2 + R_3S_2 
\end{matrix} \right).
\]
It follows from the definitions that $\mu_k(R_T^\circ)$ is
obtained from $\mu_k(R_T^\vee)$ by deleting the last $n$ columns.
This yields 
\[
\mu_k(R_T^\circ) =
\left( \begin{matrix}
S_1^tR_1S_1           \\
S_2^tR_1S_1 + R_3S_1   
\end{matrix} \right) = R_{T^*}^\circ.
\]
\end{proof}

Now we combine Corollary \ref{corrijequation} and Corollary \ref{mutation5}
and obtain the following theorem:

\begin{Thm}\label{thm_mutation}
$B(\mu_{T_k}(T))^\circ = \mu_k(B(T)^\circ)$.
\end{Thm}


\section{Examples}\label{example}


In this section we want to illustrate some of our results with examples.
We often describe modules by just indicating the multiplicities
of composition factors in the socle series.
For example for the preprojective algebra $\la$ of type
$\A_3$ we write
\[
M = {\begin{smallmatrix}&2&\\1&&3\\&2&\end{smallmatrix}}. 
\]
This means that $M$ has a socle isomorphic to the simple
labelled by $2$, the next layer of the socle series 
is isomorphic to $1 \oplus 3$, and finally the third layer is
isomorphic to $2$ again.

The examples discussed are for preprojective algebras of type
$\A_2$ and $\A_3$.
These are easy to deal with since they are representation finite
algebras.
The only other finite type case is $\A_4$ (recall that we excluded
$\A_1$), and the tame cases are $\A_5$ and $\D_4$.
Beyond that, all preprojective algebras are of wild representation
type.

\subsection{Case $\A_2$}\label{casea2}
Let $\la$ be the preprojective algebra of type $\A_2$.
There are exactly four indecomposable $\la$-modules up to isomorphism, namely
\begin{align*}
T_1 &= 
{\begin{smallmatrix}1\end{smallmatrix}}, &
T_2 &= 
{\begin{smallmatrix}2\end{smallmatrix}}, & 
T_3 &= 
{\begin{smallmatrix}1&\\&2\end{smallmatrix}}, &
T_4 &= 
{\begin{smallmatrix}&2\\1&\end{smallmatrix}}.
\end{align*}

The modules $T_3$ and $T_4$ are the indecomposable projective
$\la$-modules.
The Auslander-Reiten quiver of $\la$ looks as follows:
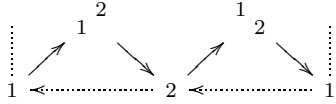
\[
\xymatrix@-1.0pc{
\ar@{.}[d]& {\begin{smallmatrix}&2\\1&\end{smallmatrix}} \ar[dr]
&&{\begin{smallmatrix}1&\\&2\end{smallmatrix}} \ar[dr] & \ar@{.}[d] \\
{\begin{smallmatrix}1\end{smallmatrix}} \ar[ur] &&
{\begin{smallmatrix}2\end{smallmatrix}} \ar[ur]\ar@{.>}[ll] &&
{\begin{smallmatrix}1\end{smallmatrix}} \ar@{.>}[ll]
}
\]
In the above picture one has to identify the two vertical dotted lines.
The dotted arrows describe the Auslander-Reiten translation.

The module $T = T_1 \oplus T_3 \oplus T_4$ is complete rigid,
and the quiver of the endomorphism algebra $\End_\la(T)$ looks as follows:
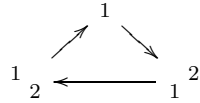
\[
\xymatrix@-1.0pc{
& {\begin{smallmatrix}1\end{smallmatrix}} \ar[dr]\\
 {\begin{smallmatrix}1&\\&2\end{smallmatrix}} \ar[ur]
&& {\begin{smallmatrix}&2\\1&\end{smallmatrix}} \ar[ll] 
}
\]

We want to illustrate Proposition \ref{homfunctor}:
The algebra $\End_\la(T)$ is isomorphic to the path algebra
of the quiver
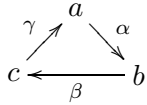
\[
\xymatrix@-1.0pc{
& a \ar[dr]^{\alpha}\\
c \ar[ur]^{\gamma}
&& b \ar[ll]^{\beta} 
}
\]
with zero relations $\beta \alpha$ and $\gamma \beta$.
There are exactly 7
indecomposable $\End_\la(T)$-modules, all of which are serial.
The Auslander-Reiten quiver of $\End_\la(T)$ looks as follows:
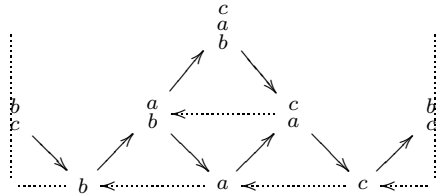
\[
\xymatrix@-1.0pc{
\ar@{.}[dd]&&& {\begin{smallmatrix}c\\a\\b\end{smallmatrix}} \ar[dr]
&&& \ar@{.}[dd]\\
{\begin{smallmatrix}&b\\&c\end{smallmatrix}} \ar[dr] &&
{\begin{smallmatrix}a\\b\end{smallmatrix}} \ar[ur] \ar[dr] &&
{\begin{smallmatrix}c\\a\end{smallmatrix}} \ar[dr] \ar@{.>}[ll]&&
{\begin{smallmatrix}b&\\c&\end{smallmatrix}}\\
&{\begin{smallmatrix}b\end{smallmatrix}} \ar[ur] \ar@{.}[l]&&
{\begin{smallmatrix}a\end{smallmatrix}} \ar[ur] \ar@{.>}[ll] &&
{\begin{smallmatrix}c\end{smallmatrix}} \ar[ur] \ar@{.>}[ll] & \ar@{.>}[l]
}
\]
Again,
the dotted arrows describe the Auslander-Reiten translation,
and the two dotted vertical lines have to be identified.

One can easily check that there are exactly four indecomposable 
$\End_\la(T)$-modules of projective dimension at most one, namely the
three indecomposable projective modules and the simple
module corresponding to the vertex $c$.
One easily checks that
\begin{align*}
F_T(T_1) &= {\begin{smallmatrix}a\\b\end{smallmatrix}}, &
F_T(T_2) &= {\begin{smallmatrix}c\end{smallmatrix}}, &
F_T(T_3) &= {\begin{smallmatrix}c\\a\\b\end{smallmatrix}}, &
F_T(T_4) &= {\begin{smallmatrix}b\\c\end{smallmatrix}}, 
\end{align*}
where $F_T$ is the functor defined in Proposition \ref{homfunctor}.

For $T$ we get
\begin{align*}
C_T &= 
\left(
\begin{matrix} 
1&1&0\\
0&1&1\\
1&1&1
\end{matrix}
\right), &
R_T &= C_T^{-t} =  
\left(
\begin{matrix} 
0&1&-1\\
-1&1&0\\
1&-1&1
\end{matrix}
\right), &
S &= S(R_T,1) = 
\left(
\begin{matrix} 
-1&0&1\\
0&1&0\\
0&0&1
\end{matrix}
\right).
\end{align*}

Besides $T$ there exists only one more basic complete rigid $\la$-module,
namely $T^* = T_2 \oplus T_3 \oplus T_4$.
The endomorphism algebras $\End_\la(T)$ and
$\End_\la(T^*)$ are isomorphic.
The two corresponding exchange sequences are
\[
0 \to {\begin{smallmatrix}2\end{smallmatrix}} \to
{\begin{smallmatrix}1&\\&2\end{smallmatrix}} \to 
{\begin{smallmatrix}1\end{smallmatrix}} \to 0
\]
and
\[ 
0 \to
{\begin{smallmatrix}1\end{smallmatrix}} \to
{\begin{smallmatrix}&2\\1&\end{smallmatrix}} \to
{\begin{smallmatrix}2\end{smallmatrix}} \to 0.
\]
For $T^*$ we get
\[
C_{T^*} = 
\left(
\begin{matrix} 
1&0&1\\
1&1&1\\
0&1&1
\end{matrix}
\right)
\text{ and }
R_{T^*} = C_{T^*}^{-t} =  
\left(
\begin{matrix} 
0&-1&1\\
1&1&-1\\
-1&0&1
\end{matrix}
\right).
\]
Now one can check that 
$SC_TS^t = C_{T^*}$ and $S^tR_TS = R_{T^*}$.

\subsection{Case $\A_3$}
Let $\la$ be the preprojective algebra of type $\A_3$.
There are exactly 12 indecomposable $\la$-modules up to isomorphism,
and all of these are rigid.

Define
$T = T_1 \oplus \cdots \oplus T_6$
where
\begin{align*}
T_1 &= 
{\begin{smallmatrix}1\end{smallmatrix}}, &
T_2 &= 
{\begin{smallmatrix}1&\\&2\end{smallmatrix}}, &
T_3 &= 
{\begin{smallmatrix}&2\\1&\end{smallmatrix}}, &
T_4 &= 
{\begin{smallmatrix}1&&\\&2&\\&&3\end{smallmatrix}}, &
T_5 &= 
{\begin{smallmatrix}&2&\\1&&3\\&2&\end{smallmatrix}}, & 
T_6 &= 
{\begin{smallmatrix}&&3\\&2&\\1&&\end{smallmatrix}}.
\end{align*}
By computing the dimension of $\End_\la(T)$ 
one can easily check that $T$ is complete rigid.
The quiver $\Gamma_T$ of the endomorphism algebra 
$\End_\la(T)$ looks as follows:
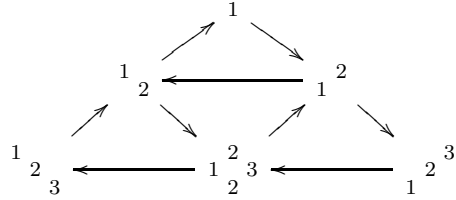
\[
\xymatrix@-1.0pc{
&& {\begin{smallmatrix}1\end{smallmatrix}} \ar[dr]\\
& {\begin{smallmatrix}1&\\&2\end{smallmatrix}} \ar[dr] \ar[ur]
&& {\begin{smallmatrix}&2\\1&\end{smallmatrix}} \ar[ll] \ar[dr]\\
{\begin{smallmatrix}1&&\\&2&\\&&3\end{smallmatrix}} \ar[ur] 
&& {\begin{smallmatrix}&2&\\1&&3\\&2&\end{smallmatrix}} \ar[ll] \ar[ur]
&& {\begin{smallmatrix}&&3\\&2&\\1&&\end{smallmatrix}} \ar[ll]
}
\]
Set 
$T^* = \mu_{T_2}(T)$.
It turns out that $T^*$ is obtained from $T$ by replacing the direct
summand $T_2$ by the module 
\[
T_2^* = 
{\begin{smallmatrix}&2&\\1&&3\end{smallmatrix}}.
\]
The two exchange sequences are
\[
0 \to {\begin{smallmatrix}&2&\\1&&3\end{smallmatrix}} \to 
{\begin{smallmatrix}&2\\1&\end{smallmatrix}}  \oplus 
{\begin{smallmatrix}1&&\\&2&\\&&3\end{smallmatrix}} \to 
{\begin{smallmatrix}1&\\&2\end{smallmatrix}} \to 0
\]
and
\[
0 \to {\begin{smallmatrix}1&\\&2\end{smallmatrix}} \to 
{\begin{smallmatrix}1\end{smallmatrix}} \oplus 
{\begin{smallmatrix}&2&\\1&&3\\&2&\end{smallmatrix}} \to 
{\begin{smallmatrix}&2&\\1&&3\end{smallmatrix}}  \to 0.
\]
The quiver $\Gamma_{T^*}$ of $\End_\la(T^*)$ looks as follows:
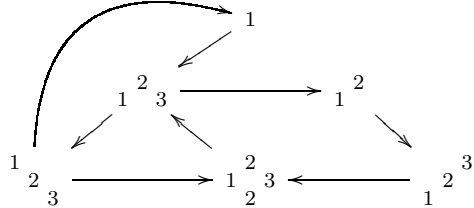
\[
\xymatrix@-1.0pc{
&& {\begin{smallmatrix}1\end{smallmatrix}} \ar[dl]\\
& {\begin{smallmatrix}&2&\\1&&3\end{smallmatrix}} \ar[rr] \ar[dl]
&& {\begin{smallmatrix}&2\\1&\end{smallmatrix}} \ar[dr]\\
{\begin{smallmatrix}1&&\\&2&\\&&3\end{smallmatrix}} \ar[rr] \ar@/^3pc/[uurr]
&& {\begin{smallmatrix}&2&\\1&&3\\&2&\end{smallmatrix}} \ar[ul]
&& {\begin{smallmatrix}&&3\\&2&\\1&&\end{smallmatrix}} \ar[ll]
}
\]

For $T$ an easy calculation yields 
\[
C_T = 
\left(
\begin{matrix} 
1&1&0&1&0&0\\
0&1&1&1&1&0\\
1&1&1&1&1&0\\
0&0&0&1&1&1\\
0&1&1&1&2&1\\
1&1&1&1&1&1
\end{matrix}
\right)
\text{ and }
R_T = C_T^{-t} =  
\left(
\begin{matrix} 
0&1&-1&0&0&0\\
-1&0&1&1&-1&0\\
1&-1&0&0&1&-1\\
0&-1&0&1&0&0\\
0&1&-1&-1&1&0\\
0&0&1&0&-1&1
\end{matrix}
\right),
\]
and for $T^*$ we get
\[
C_{T^*} = 
\left(
\begin{matrix} 
1&0&0&1&0&0\\
1&1&0&1&1&1\\
1&1&1&1&1&0\\
0&1&0&1&1&1\\
0&1&1&1&2&1\\
1&1&1&1&1&1
\end{matrix}
\right)
\text{ and }
R_{T^*} = C_{T^*}^{-t} =  
\left(
\begin{matrix} 
0&-1&0&1&0&0\\
1&0&-1&-1&1&0\\
0&1&0&0&0&-1\\
-1&1&0&1&-1&0\\
0&-1&0&0&1&0\\
0&0&1&0&-1&1
\end{matrix}
\right).
\]
Furthermore, we have
\[
S = S(R_T,2) =
\left(
\begin{matrix} 
1&&&&&\\
1&-1&&&1&\\
&&1&&&\\
&&&1&&\\
&&&&1&\\
&&&&&1
\end{matrix}
\right).
\]
In the above matrix we only displayed the non-zero entries of $S$.
Now one can check that
$SC_TS^t = C_{T^*}$ and $S^tR_TS = R_{T^*}$.

For the case $\A_3$ there are exactly 14 basic complete rigid 
$\la$-modules up to isomorphism.

\subsection{Case $\A_4$}
Let $\la$ be the preprojective algebra of type $\A_4$.
There are exactly 40 indecomposable $\la$-modules up to isomorphism,
and all of these are rigid.
The number of isomorphism classes of basic complete rigid $\la$-modules
is 672.
For more details we refer to \cite{GLS}.


\section{Relation with cluster algebras}\label{cluster}


From now on let $\field = \C$ be the field of complex numbers.

\subsection{}
We are going to outline the construction of the map
$\varphi : M\mapsto\varphi_M$ introduced in~\ref{lift}.
As before let $\Lam_\be$ be the affine variety of 
$\Lam$-modules with dimension vector $\be \in \N^n$.
Let $\MB = \oplus_{\be\in\N^n} \MB_\be$ be the
algebra of $G_\be$-invariant constructible functions 
introduced by Lusztig as a geometric model for $U(\n)$.
Here $\MB_\be$ is the vector space of functions from $\Lam_\be$
to $\C$ spanned by certain functions $d_\ib$ defined
as follows.
For $M\in\Lam_\be$, let $\Phi_{\ib,M}$ be the variety of 
composition series of $M$ of type $\ib=(i_1,\ldots,i_k)$, 
that is of flags of submodules 
\[
M=M_0 \supset M_1 \supset \cdots \supset M_k = {0}
\]
with $M_{j-1}/M_j$ isomorphic to the simple module 
$S_{i_j}$ for $j=1,\ldots ,k$.
Then 
\begin{equation}
d_\ib(M):=\chi(\Phi_{\ib,M}),
\end{equation}
where $\chi$ denotes the Euler characteristic.

Let $\MB^*= \oplus_{\be\in\N^n} \MB_\be^*$ be the graded
dual of $\MB$, and let $\delta_M\in\MB^*$ be the linear
form which maps a constructible 
function $f\in\MB$ to its evaluation $f(M)$ at $M$.
It is completely determined by the numbers
$\delta_M(d_\ib) = \chi(\Phi_{\ib,M})$ where $\ib$
varies over the possible composition types. 
Under the isomorphism $\MB^* \cong U(\n)^* \cong \C[N]$,
$\delta_M$ gets identified with a regular function 
$\varphi_M\in\C[N]$. 

To describe $\varphi_M$ explicitly, we introduce the
one-parameter subgroups
\[
x_i(t) = \exp(te_i), \quad (t\in \C,\ i\in Q_0),
\]
of $N$ associated with the Chevalley generators $e_i$
of $\n$.

\begin{Lem}
For any sequence $\ib=(i_1,\ldots,i_k)$ of elements of $Q_0$,
we have
\[
\varphi_M(x_{i_1}(t_1)\cdots x_{i_k}(t_k))
=
\sum_{\ab=(a_1,\ldots,a_k)\in \N^k}
\chi(\Phi_{\ib^\ab,M}) 
\frac{t_1^{a_1}}{a_1!} \cdots \frac{t_k^{a_k}}{a_k!}\ ,
\qquad (t_1,\ldots,t_k \in \C),
\]
where we use the short-hand notation 
$\ib^\ab:=(\underbrace{i_1,\ldots,i_1}_{a_1},\ldots,
\underbrace{i_k,\ldots,i_k}_{a_k})$.
\end{Lem}

The proof follows easily from the classical description
of the duality between $U(\n)$ and $\C[N]$.
Note that if $\ib$ is a reduced word for the longest
element of the Weyl group of $\g$ then the set 
$\{ x_{i_1}(t_1)\cdots x_{i_k}(t_k)\mid t_1,\ldots,t_k \in \C\}$ 
is dense in $N$, so the above formula completely determines the polynomial
function $\varphi_M$.

It is shown in \cite[Lemma 7.3]{GLS} that the functions $\varphi_M$
are multiplicative, in the sense that
\begin{equation}\label{multiplicativity}
\varphi_M\, \varphi_N = \varphi_{M\oplus N},\qquad (M,N\in \md(\Lam)).
\end{equation}

\subsection{}\label{semcan}
Let $Z$ be an irreducible component of the variety $\Lam_\be$.
The map $\varphi$ being constructible, there exists
a dense open subset $O_Z$ of $Z$ such that for all $M,N\in O_Z$
we have $\varphi_M = \varphi_N$.
A point $M\in O_Z$ is called a {\it generic point} of $Z$.
Put $\rho_Z=\varphi_M$ for a generic point $M$ of $Z$.
Then the collection $\{\rho_Z\}$ where $Z$ runs over all
irreducible components of all varieties $\Lam_\be$ is dual
to the semicanonical basis of $U(\n)$. We call it the 
{\it dual semicanonical basis} of $\C[N]$ and denote it by
$\SB^*$.
 
If $M\in\Lam_\be$ is a rigid $\Lam$-module in the
irreducible component $Z$, its $G_\be$-orbit is open 
so $M$ is generic and $\varphi_M=\rho_Z$ belongs to $\SB^*$.

\subsection{}
Recall the setting of \ref{clusterCN}.
From the complete rigid $\Lam$-module
$T_Q=T_1\oplus\cdots\oplus T_r$ of 
Theorem~\ref{startmoduletheorem}
we get the initial seed $(\xb(T_Q), B(T_Q)^\circ)$ where
\[ 
\xb(T_Q)=(x_1(T_Q),\ldots,x_r(T_Q))
:=(\varphi_{T_1},\ldots,\varphi_{T_r}).
\] 
By \cite{GLS2}, it coincides with one
of the initial seeds of \cite{BFZ} for the upper cluster
algebra structure on $\C[N]$. 
Here we assume as before that $T_{r-n+1}, \ldots , T_r$
are the $n$ indecomposable projective $\Lam$-modules,
so that $x_1(T_Q),\ldots,x_{r-n}(T_Q)$ are the exchangeable 
cluster variables of $\xb(T_Q)$ and the remaining ones generate
the ring of coefficients.

Recall the definition of the mutation of a seed 
$(\xb,\tB)$ in direction $k\in [1,r-n]$.
This is the new seed $(\xb',\tB')$, where 
$\tB'=\mu_k(\tB)$ is given by \ref{matmutation},
and $\xb'$ is obtained from $\xb=(x_1,\ldots,x_r)$ by replacing
$x_k$ by
\begin{equation}\label{exchrel}
x_k'= \frac{\prod_{b_{ik}>0} x_i^{b_{ik}}
+ \prod_{b_{ik}<0} x_i^{-b_{ik}}}{x_k}.
\end{equation}  
Here the exponents $b_{ik}$ are the entries of the matrix $\tB$.

The mutation class $\CB$ of the seed $(\xb(T_Q), B(T_Q)^\circ)$ is 
defined to be the set of all seeds $(\xb,\tB)$ which 
can be obtained from $(\xb(T_Q), B(T_Q)^\circ)$ by a sequence of mutations.

\subsection{} 
We will need the following result of \cite{GLS3}.

\begin{Thm}\label{multform}
Let $M$ and $N$ be $\Lam$-modules such that
$\dm \Ext_\la^1(M,N) = 1$,
and let 
\[
0 \to M \to X \to N \to 0
\text{\;\;\; and \;\;\;}
0 \to N \to Y \to M \to 0
\]
be non-split short exact sequences.
Then
$\varphi_{M} \cdot \varphi_{N} = \varphi_X + \varphi_Y$.
\end{Thm}

\subsection{} 
Now everything is ready for the proof of 
Theorem~\ref{clustermonom}.
Let $R=R_1\oplus\cdots\oplus R_r$ be a vertex of
$\T_\Lam^\circ$. Thus $R$ is obtained from $T_Q$ by means
of a finite number of mutations, say $\ell$.
We want to prove that $(\xb(R),B(R)^\circ)$ is a seed in $\CB$.  
We argue by induction on $\ell$.
Put $\xb(R)=(\varphi_{R_1},\ldots,\varphi_{R_r})$.
If $\ell=0$ then $R=T_Q$ and we know already that 
$(\xb(R),B(R)^\circ)$ is an initial seed of $\CB$.
Otherwise we have $R=\mu_k(S)$ for some vertex $S$
of $\T_\Lam^\circ$ and some $k$, and by induction
we can assume that $(\xb(S),B(S)^\circ)$ is a seed in $\CB$. 
By Corollary~\ref{corquivershape}, we know that
$\dm \Ext_\la^1(S_k,R_k) = 1$, so we can apply
Theorem~\ref{multform}.
Let 
\[
0 \to R_k \to X \to S_k \to 0
\text{\;\;\; and \;\;\;}
0 \to S_k \to Y \to R_k \to 0
\]
be non-split short exact sequences.
Then $\varphi_{R_k} \cdot \varphi_{S_k} = \varphi_X + \varphi_Y$.
By Corollary~\ref{corrijequation} and the remarks 
which follow it, we have 
\[
X = \bigoplus_{b_{ik}<0}S_i^{-b_{ik}},
\qquad
Y = \bigoplus_{b_{ik}>0}S_i^{b_{ik}},
\]
where the $b_{ik}$ here are the entries of $B(S)^\circ$.
Using (\ref{multiplicativity}), it follows that
\[
\varphi_X = \prod_{b_{ik}<0} \varphi_{S_i}^{-b_{ik}},
\qquad
\varphi_Y = \prod_{b_{ik}>0} \varphi_{S_i}^{b_{ik}}.
\]
Hence, comparing with (\ref{exchrel}) and taking 
into account Theorem~\ref{thm_mutation}, we see that 
$(\xb(R),B(R)^\circ)$ is obtained from $(\xb(S),B(S)^\circ)$
by a seed mutation in direction $k$.
This shows that the map $R \mapsto (\xb(R),B(R)^\circ)$ gives a covering
of graphs from $\T_\Lam^\circ$ to the exchange graph $\G$
of $\C[N]$.

Now if $R$ and $R'$ are such that $\xb(R)=\xb(R')$, then
$\varphi_R = \varphi_{R'}$. In particular, $R$ and $R'$
have the same dimension vector $\be$.
Since $R$ and $R'$ are both generic $\varphi_R = \varphi_{R'}$
belongs to $\SB^*$ and $R$ and $R'$ have to be in the same
irreducible component of $\Lam_\be$. Therefore $R$ and $R'$
are isomorphic. Hence the covering of graphs is in fact
an isomorphism.
(Two seeds which can be obtained from each other by reordering of
the entries of the clusters and a corresponding reordering of the
columns and rows of the exchange matrices are considered to be identical.)

Finally, using (\ref{multiplicativity}), we get that all cluster
monomials are of the form $\varphi_M$ for a rigid module 
$M$ (not necessarily basic or maximal), hence by \ref{semcan}
they belong to $\SB^*$.

This finishes the proof of Theorem~\ref{clustermonom}.

\bigskip
{\parindent0cm \bf Acknowledgements.}\,
We like to thank William Crawley-Boevey, Osamu Iyama, Bernt Jensen, 
Bernhard Keller and Idun Reiten
for interesting and helpful discussions.
We thank Iyama for explaining us the results in Section \ref{iyamaresults} 
(which were at that time not yet available in a written form) during a 
conference \cite{I2} in August 2004.


\end{document}